\documentclass[12pt,times]{amsart}
\usepackage{graphicx}
\usepackage{amssymb}
\numberwithin{equation}{section}

\date{}
\newcommand{\calA}{\mathcal{A}}
\newcommand{\calB}{\mathcal{B}}
\newcommand{\calC}{\mathcal{C}}
\newcommand{\calD}{\mathcal{D}}

\newcommand{\calH}{\mathcal{H}}

\newcommand{\calL}{\mathcal{L}}

\newcommand{\calN}{\mathcal{N}}
\newcommand{\calO}{\mathcal{O}}
\newcommand{\calP}{\mathcal{P}}
\newcommand{\calQ}{\mathcal{Q}}

\newcommand{\calS}{\mathcal{S}}
\newcommand{\calI}{\mathcal{I}}
\newcommand{\calT}{\mathcal{T}}

\newcommand{\bbF}{\mathbb{F}}
\newcommand{\bbC}{\mathbb{C}}

\newcommand{\bbP}{\mathbb{P}}

\newcommand{\bbR}{\mathbb{R}}

\newcommand{\bbZ}{\mathbb{Z}}

\newcommand{\bfF}{\mathbf{F}}

\def\MW{\textup{MW}}
\def\SL{\textup{SL}}
\def\Aut{\textup{Aut}}

\def\PG{\textup{PG}}
\def\Sym{\textup{Sym}}
\def\AG{\textup{AG}}

\def\PGL{\textup{PGL}}
\def\Ceva{\textup{Ceva}}
\def\Sp{\textup{Sp}}
\def\Kum{\textup{Kum}}

\title{ABSTRACT CONFIGURATIONS IN ALGEBRAIC GEOMETRY}
\author{I. Dolgachev}
\address{Department of Mathematics, University of Michigan, Ann Arbor, MI 48109,USA}
\email{idolga@umich.edu}
\dedicatory{To the memory of Andrei Tyurin}
\begin{document}

\begin{abstract}
An abstract  $(v_k,b_r)$-configuration is a pair of finite sets of cardinalities $v$ and $b$ with a relation on the product of the sets  such that each element of the first set is related to the same number $k$ of elements from the second set and, conversely, each element of the second set is related to the same number $r$ of elements in the first set. An example of an abstract configuration is a finite geometry. In this paper we discuss some examples of abstract configurations and, in particular finite geometries, which one encounters in algebraic geometry.
\end{abstract}

\maketitle

CONTENTS
\begin{itemize}
\item[1.] Introduction
\item[2.] Configurations, designs and finite geometries
\item[3.] Configurations in algebraic geometry 
\item[4.] Modular configurations
\item[5.] The Ceva configurations
\item[6.] $v_3$-configurations
\item[7.] The Reye $(12_4,16_3)$-configuration
\item[8.] $v_{v-1}$-configurations
\item[9.] The Cremona-Richmond $15_3$-configuration
\item[10.] The Kummer configurations
\item[11.] A symmetric realization of $\bbP^2(\bbF_q)$
\end{itemize}

\section{Introduction} In this paper we discuss some examples of abstract configurations and, in particular finite geometries, which one encounters in algebraic geometry. The Fano  Conference makes it very appropriate because of the known contribution of Gino Fano to  finite geometry. In his first published paper \cite{Fano1}  he gave a first  synthetical definition of the projective plane over an arbitrary field. In \cite{Fano1} introduces his famous Fano's Postulate. It asserts that the diagonals of a complete quadrangle do not intersect at one point. Since this does not hold for the projective plane over a finite field of two elements (the Fano plane), the postulate  allows one to exclude the case of characteristic 2. Almost forty years later he returned to the subject of finite projective planes in \cite{Fano2}.

An abstract  $(v_k,b_r)$- configuration is a pair of finite sets of cardinalities $v$ and $b$ with a relation on the product of the sets  such that each element of the first set is related to the same number $k$ of elements from the second set and, conversely, each element of the second set is related to the same number $r$ of elements in the first set. A geometric realization of an abstract configuration is a realization of  each set as a set of linear (or affine, or projective) subspaces of certain dimension with a certain incidence relation. Any abstract configuration admits a geometric realization in a projective space of sufficient large dimension over an infinite field. The existence of a  geometric realization in a given space over a given field, for example, by real points and lines in the plane, is a very difficult problem.  A regular configuration is an abstract configuration which admits a group of automorphisms which acts transitively on the sets $\calA$ and $\calB$.  Many regular configurations arise in algebraic geometry, where the sets $\calA$ and $\calB$ are represented by subvarieties of an algebraic variety with an appropriate incidence relation. The most notorious example is the Kummer configuration $(16_6,16_6)$ of points and tropes of a Kummer surface, or the Hesse configuration $(9_3,12_4)$ arising from  inflection points of a nonsingular plane cubic. In this paper we discuss other well-known or less known configurations which admit an ``interesting'' realizations in algebraic geometry.  I am not attempting to define what does the latter mean.  I leave to the reader to decide whether the realization is an interesting realization or not. 
	
 I am thankful to J. Keum for some critical comments on the paper and to D. Higman and V. Tonchev for helping with references to the theory of designs.

\section{Configurations, designs  and finite geometries}
 \subsection{}An abstract  \emph{configuration}  is a triple $\{\calA,\calB,R\}$, where $\calA, \calB$ are non-empty finite sets and $R\subset \calA\times\calB$ is a relation such that the cardinality of the set
$$R(x) = \{B\in \calB:(x,B)\in R\}$$
(resp, the set
$$R(B) = \{x\in \calA:(x,B)\in R\})$$
does not depend on  $x\in \calA$ (resp.  $B\in \calB$). Elements of $\calA$ are called \emph{points}, elements of $\calB$ are called \emph{blocks}.   If 
$$v= \#\calA, \quad b = \#\calB, \quad k =\#R(x),\quad  r = \#R(B),$$
 then a configuration is said to be an $(v_k,b_r)$-configuration.  A {\it symmetric configuration} is a configuration with $ \#\calA = \#\calB$, and hence $k = r$. It is said to be an $v_k$-configuration. 

We will assume that for any two distinct elements  $x,y\in\calA$ (resp.
$B,B'\in\calB)$, we have $R(x)\ne R(y)$ (resp. $R(B)\ne R(B')$).  This, for example, excludes symmetric configurations $v_k$ with $v =k$. Under this assumption we can identify the set $\calB$ with a set of subsets (blocks) of the set $\calA$ of cardinality $r$  such that each element of $\calA$ belongs to exactly $k$ blocks.

Note that we consider only a special case of the notion of an abstract  configuration studied in combinatorics. Ours is  a {\it tactical configuration}.

\subsection{The Levi graph} An abstract $(v_k,b_r)$-configuration can be uniqu\\ely represented by its {\it Levi graph} (see \cite{Cox}). In this graph elements of $\calA$ are represented by black vertices and elements of $\calB$ by  white vertices. An edge joins a  black vertex and a white vertex if and only if the corresponding elements are in the relation. The Levi graph of a $v_k$-configuration is $k$-valent, i.e. each vertex is incident with exactly $k$ edges. A $k$-valent graph with $2v$-vertices is the Levi graph of a $v_k$-configuration if the set of vertices can be partitioned into two sets of cardinality $v$ such that two vertices belonging to one set are not connected by an edge.  

\subsection{Direct sums, complements}\label{oper} One defines naturally the direct sum of configurations of type $(v_k,b_r)$ and type 
$(v'_{k'},b'_{r'})$. This is a configuration of type $\bigl((v+v')_{k},(b+b')_{r}\bigr)$. A configuration is called {\it connected} if it is not equal to a direct sum of configurations. Or, equivalently,   for any $x,y\in \calA\cup \calB$ there exist $z_1,\ldots,z_k\in \calA\cup \calB$ such that $(xRz_1), (z_1Rz_2,\ldots,(z_nR,y)$, where for any $x,y\in \calA\cup\calB$ we write $xRy$ if either $(x,y)\in R$ or $(y,x)\in R$.

For any symmetric $v_k$-configuration, one defines the {\it complementary} $(v_{v-k})$-configuration whose blocks are the complementary sets of blocks of the first configuration. 

\subsection{Symmetry} Let $\Sym(\calC)$ denote the group of symmetries of a  $(v_k,b_r)$-configuration $\calC$  defined as a subgroup of bijections $g$ of the set $\calA\cup \calB$ such that for any $p\in \calA, q\in \calB,$ we have $g(p)\in R(g(q))$ if and only if $p\in R(q)$ and $g(q)\in R(g(p))$ if and only if $q\in R(p)$.  Note that, if $k\ne r$, the subsets $\calA$ and $\calB$ are necessarily invariant under any symmetry (since $g(R(p)) = R(g(p))$). However, if the configuration is symmetric, $\Sym(\calC)$ may contain an element such that $g(\calA) = \calB$. We will call it a {\it switch}. A switch of order 2 is called a  {\it polarity}, or {\it duality}, see (\cite{Dem}). The group of symmetries of the configuration which preserve each set $\calA$ and $\calB$ is of index $\ge 2$ in the full group of symmetries. We will call it the {\it proper group of symmetries} of the configuration. 

If the configuration is connected, then each symmetry  either preserves the sets $\calA$ and $\calB$, or a switch. Of course, the notion of a symmetry of a configuration is a special case  of the notion of an  {\it isomorphism of configurations}, or even a {\it morphism of configurations}.

A configuration is called {\it $s$-regular} if its symmetry group acts transitively on the set of $s$-arcs  but not on the set of $(s+1)$-arcs of its Levi graph. Recall that an $s$-arc of a graph is an ordered sequence of oriented $s$-edges consecutively incident.

\subsection{Designs} A \emph{block-scheme} or a \emph{design} is a $(v_k,b_r)$-configuration such that, for any distinct $x,y\in \calA$, the cardinality $\lambda$ of the set $R(x)\cap R(y)$
does not depend on $x,y$ and is non-zero. The number
$\lambda$ is called the \emph{type} of a $(v_k,b_r)$-design. A standard argument shows that
\begin{equation}\label{design}
vk  = br,\quad k(r-1) =\lambda (v-1).
\end{equation}
In particular, $\lambda$ is determined by $v,k,b,r$. 

A symmetric $v_k$-design of type $\lambda$ has the additional property that, for  any $B_1,B_2\in \calB$, $\# R(B_1)\cap R(B_2) = \lambda$.

 According to a theorem of Bruck-Chowla-Ryser (see \cite{Hall}) $k-\lambda$ is always  a square if $v$ is even and the quadratic form $(k-\lambda)x^2+(-1)^{\frac{v-1}{2}}\lambda y^2-z^2$ represents zero if $v$ is odd.

\subsection{Projective geometries} (see \cite{Hirsch2}) Let $\bbP^{n}(\bbF_q)$ be the set of points of projective space of  dimension $n$ over a finite field of $q = p^r$ elements. Let $G_{r,n}$ be the Grassmannian of $r$-dimensional subspaces in $\bbP^n$. It is known that 
\[
\# G_{r,n}(\bbF_q) = m(r,n;q), \]
where
\begin{equation}\label{gras}
m(r,n;q) = \Bigl[\begin{matrix}n+1\\
r+1\end{matrix}\Bigr]_q: = \frac{(1-q^{n+1})\cdots(1-q^{n-r+1})}{(1-q)\cdots (1-q^{r+1})}.
\end{equation}
Let $\calA = G_{r,n}(\bbF_q), \calB = G_{s,n}(\bbF_q)$, $r < s$.  We take for $R$ the incidence relation $R = \{(L,M)\in \calA\times \calB:L\subset M\}$. 
Since each $r$-subspace is contained in $m(s-r-1,n-r-1;q)$ subspaces of dimension $s$, and each subspace of dimension $s$ contains $m(r,s;q)$ subspaces of dimension $r$, we obtain a $(m(r,n;q)_{m(s-r-1,n-r-1;q)},m(s,n;q)_{m(r,s;q)})$-configuration. 

We say that it is a \emph{ projective geometry configuration} and denote it by $\PG(n,r,s;q)$.  It is a design if and only if $r=0$. The type $\lambda$ is equal to $m(s-2,n-2;q)$.

Note that the projective duality gives an isomorphism of configurations
\begin{equation}\label{duality}
\PG(n,r,s;q) \cong \PG(n,n-s-1,n-r-1;q).
\end{equation}
In the case $s = n-r-1$, we get a symmetric $m(r,n;q)_{m(s-r-1,n-r-1;q)}$-configuration. 

The smallest design obtained in this way is $\PG(2,0,1;2)$. It is called the \emph{Fano plane}. 

\begin{figure}[h]
\begin{center}
\includegraphics[width=2.5in]{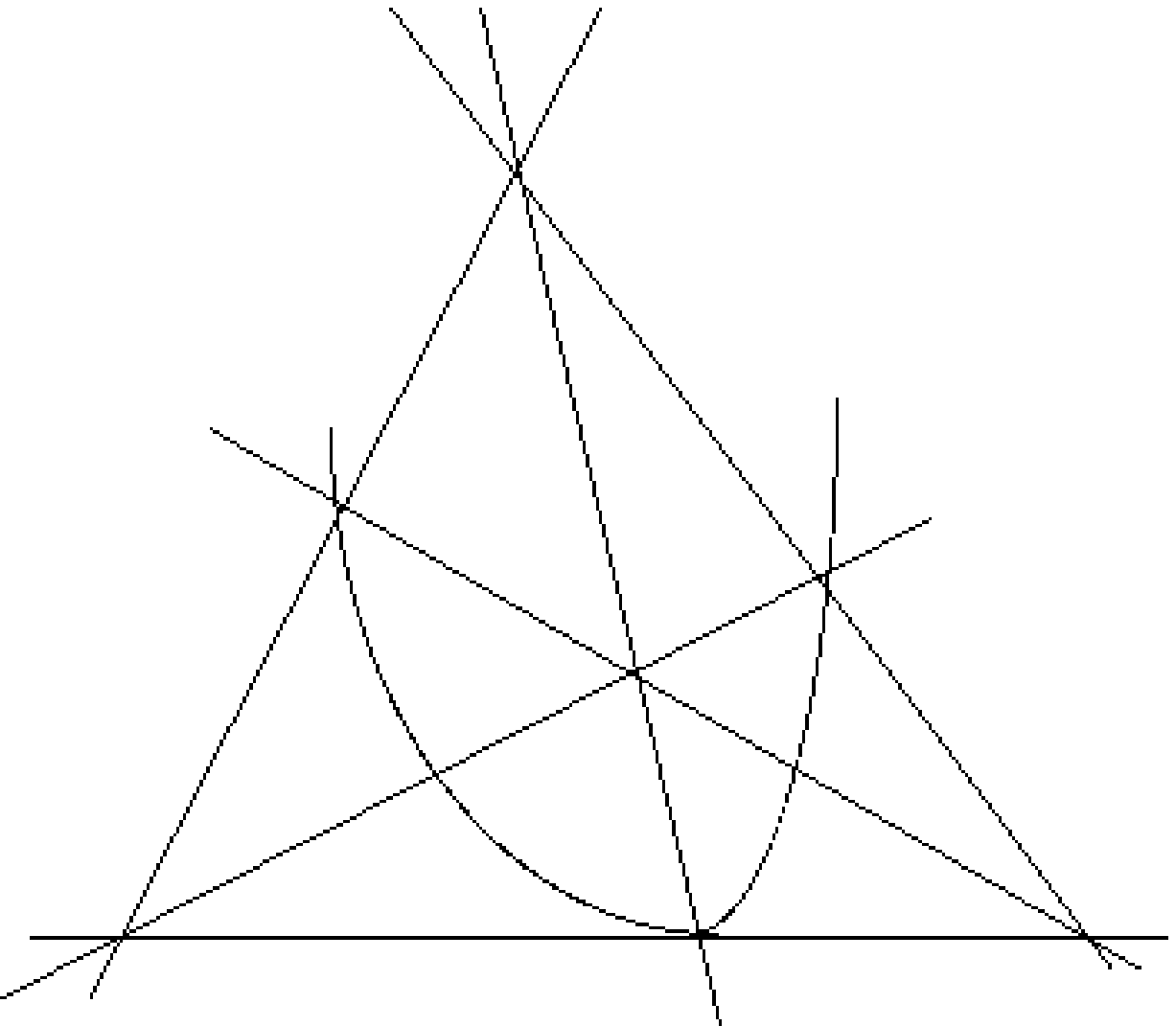}
\caption{ }
\label{}
\end{center}
\end{figure}

Its Levi graph is a  $4$-regular graph (its 4-arcs correspond to sequences point-line-point-line-point or line-point-line-point-line).

\begin{figure}[h]
\begin{center}
\includegraphics[width=2in]{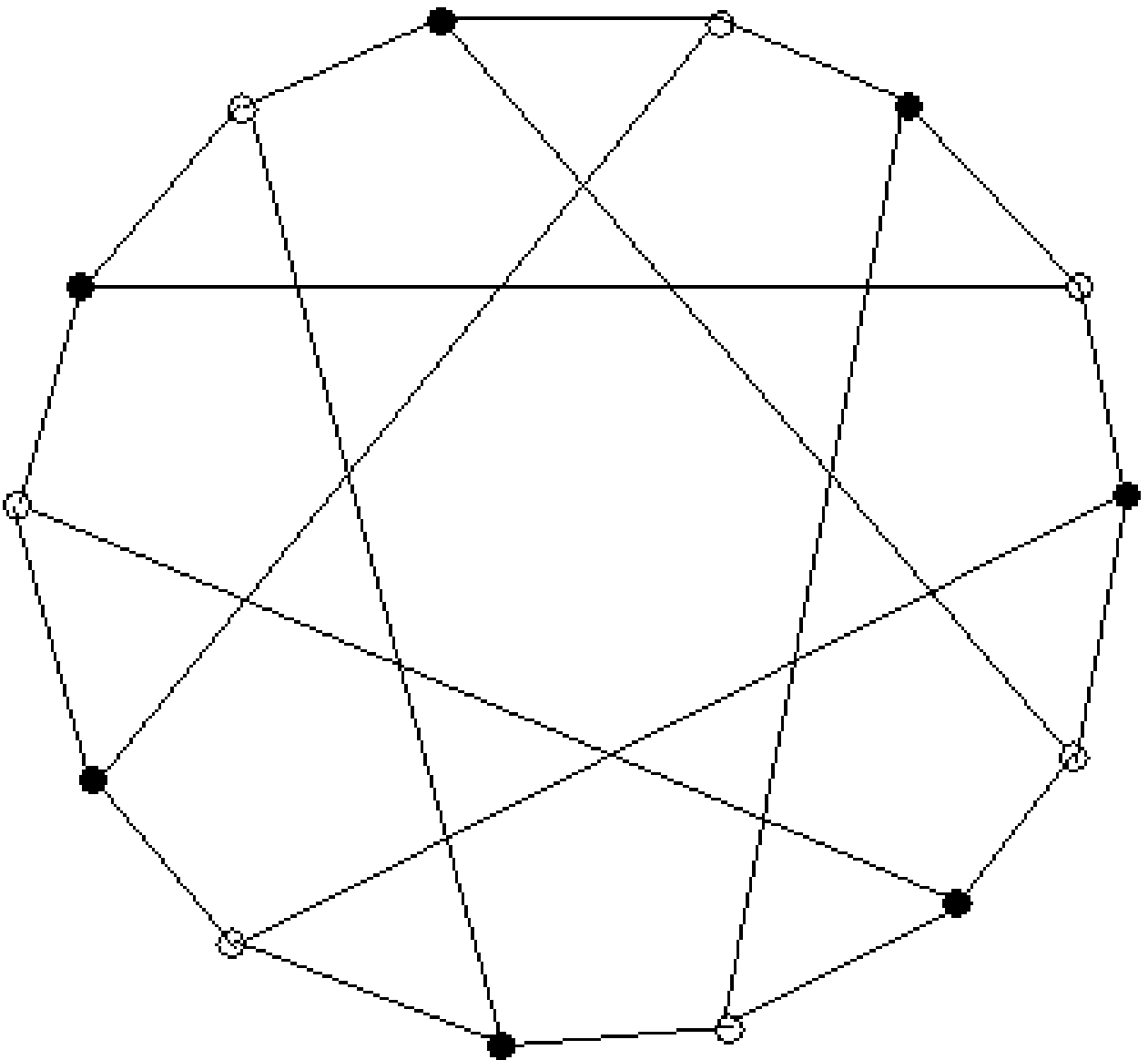}
\caption{ }
\label{fanoplane}
\end{center}
\end{figure}

\subsection{An abstract projective geometry configurations} One defines an abstract version of a configuration $\PG(n,0,1;q)$ as a $(v_k,b_r)$-design with $\lambda  = 1$ and $b\ge 3$.  We call the elements of $\calB$ lines and,  for any $x,y\in \calP$ denote the unique line in  $R(x,y)$ by $x+y$. If $x\in R(B)$, we say that the line $B$ contains the point $x$. A design of type $\lambda = 1$ is called a {\it line-point} abstract configuration (see \cite{Levi}).

A \emph{projective geometry} is a line-point configuration satisfying the following additional  condition:

(*) If $x\in B_1\cap B_2$ and $x',x''\in B_1\setminus\{x\}, y',y''\in B_2\setminus\{x\}$, then the lines $x'+y'$ and $x''+y''$ have  a common point.

\begin{figure}[h]
\begin{center}
\includegraphics[width=2in]{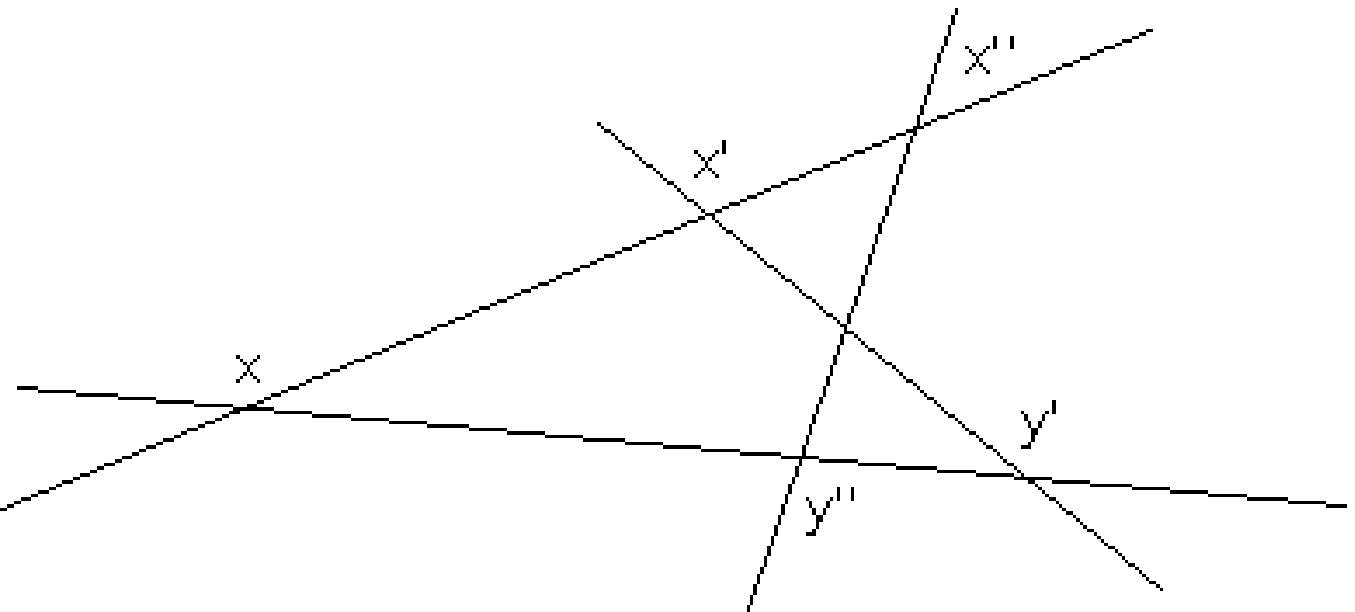}
\caption{ }
\label{ }
\end{center}
\end{figure}
One defines a subspace of a projective geometry as a subset $L$ of $\calP$ such that, 
for any distinct $x,y\in L$, any point $z\in x+y$ belongs to $L$. This allows one to define the \emph{dimension} of a projective geometry. A theorem due to D. Hilbert says that  that all projective geometries of dimension $> 2$ are isomorphic to a configuration $\PG(n,0,1;q)$ or its affine analog $\AG(n,0,1,q)$ (see \cite{Dem}, p.26). A projective plane is  
 isomorphic to $\PG(2,0,1;q)$ or its affine analog if and only if the following Desargues axiom is satisfied.

(**) Let $(x,y,z)$ and $(x',y',z')$ be two triangles without common vertices (i.e. ordered sets of three points with no common points). Assume that the 3 pairs of lines $(x+y, x'+y')$, $(y+z,y'+z')$ and $(x+z,x'+z')$ intersect at three collinear points. Then the lines $x+x',y+y', z+z'$ have a common point.

\begin{figure}[h]
\begin{center}
\includegraphics[width=2in]{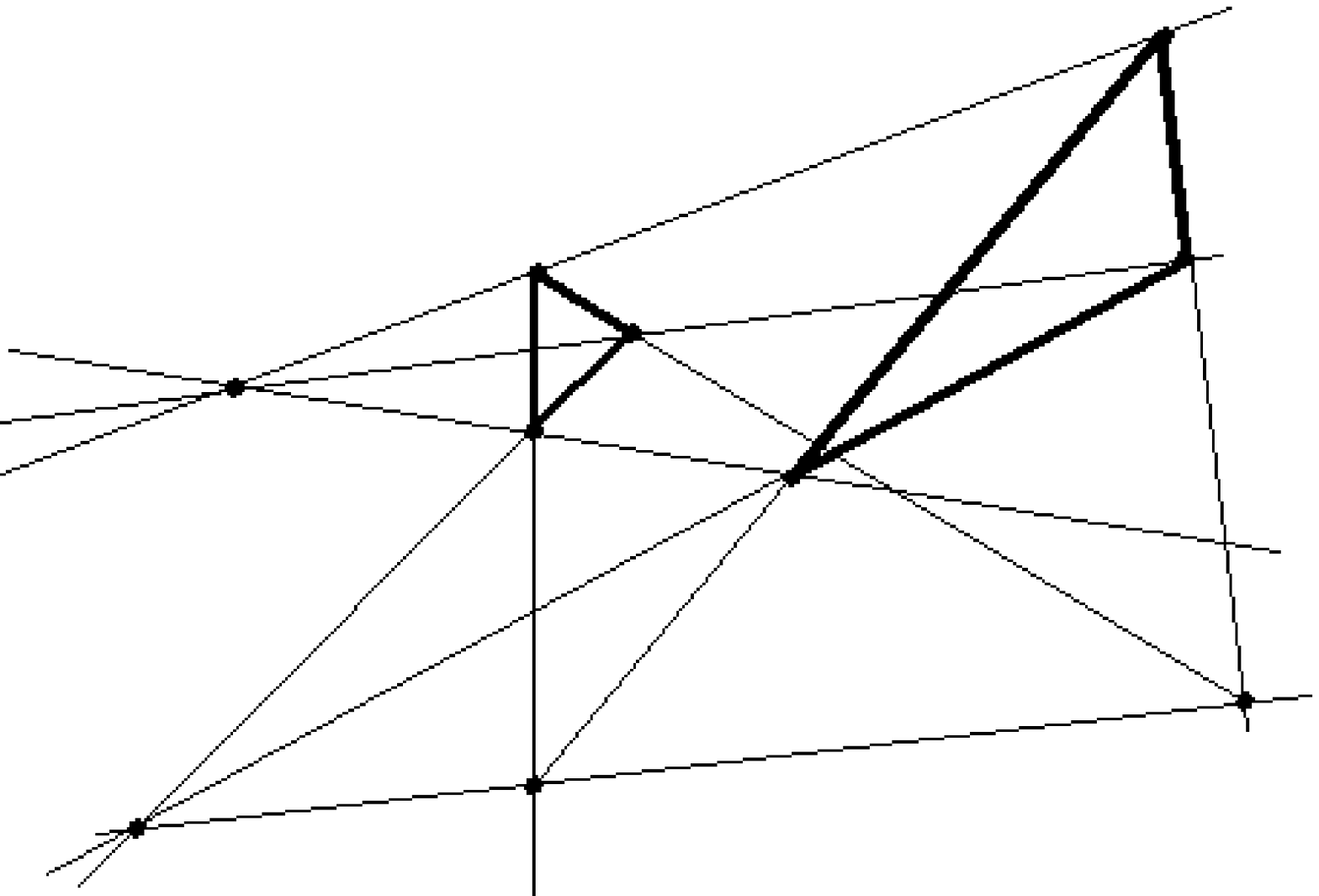}
\caption{ }
\label{ }
\end{center}
\end{figure}

Any $v_k$-design of type $\lambda = 1$ is defined by an abstract projective plane (see \cite{Hugh}, Theorem 2.21). It follows from equation \eqref{design} that $v = (k-1)^2+(k-1)^2+1$.

\subsection{}\label{desargconf} Note that, when the condition of the Desargues axiom are satisfied we have 10 lines and 10 points forming a symmetric configuration $10_3$. It is called the {\it Desargues configuration}.  An easiest  way to construct such configuration is as follows. Take a basis $e_1,\ldots,e_4$ of a vector space $V$ over a field $k$ and add a vector $e_5 = e_1+e_2+e_3+e_4$. Let $\pi_{ijk}$ be the planes  in $\bbP(V)$ spanned by $e_i, e_j, e_k$. Let $L_{ij}$ be the line spanned  by the vectors $e_i,e_j$.
A general plane $\pi$ in $\bbP(V)$ intersects the union of the 10 planes $\pi_{ijk}$ at the union of 10 lines $l_{ijk}$. Each line contains three points $p_{ij},p_{ik},p_{ik}$. They are the points of intersection of $\pi$ with the lines $L_{ij}$, $L_{ik}, L_{jk}$. Each point $p_{ij}$ is contained  in three lines $l_{ijm},l_{ijn}, l_{ijt}$, where $\{m,n,t\}$ is the subset of $\{1,\ldots,5\}$ complementary to $\{i,j\}$. Here is the Levi graph of the Desargues configuration:
\begin{figure}[h]
\includegraphics[width=2in]{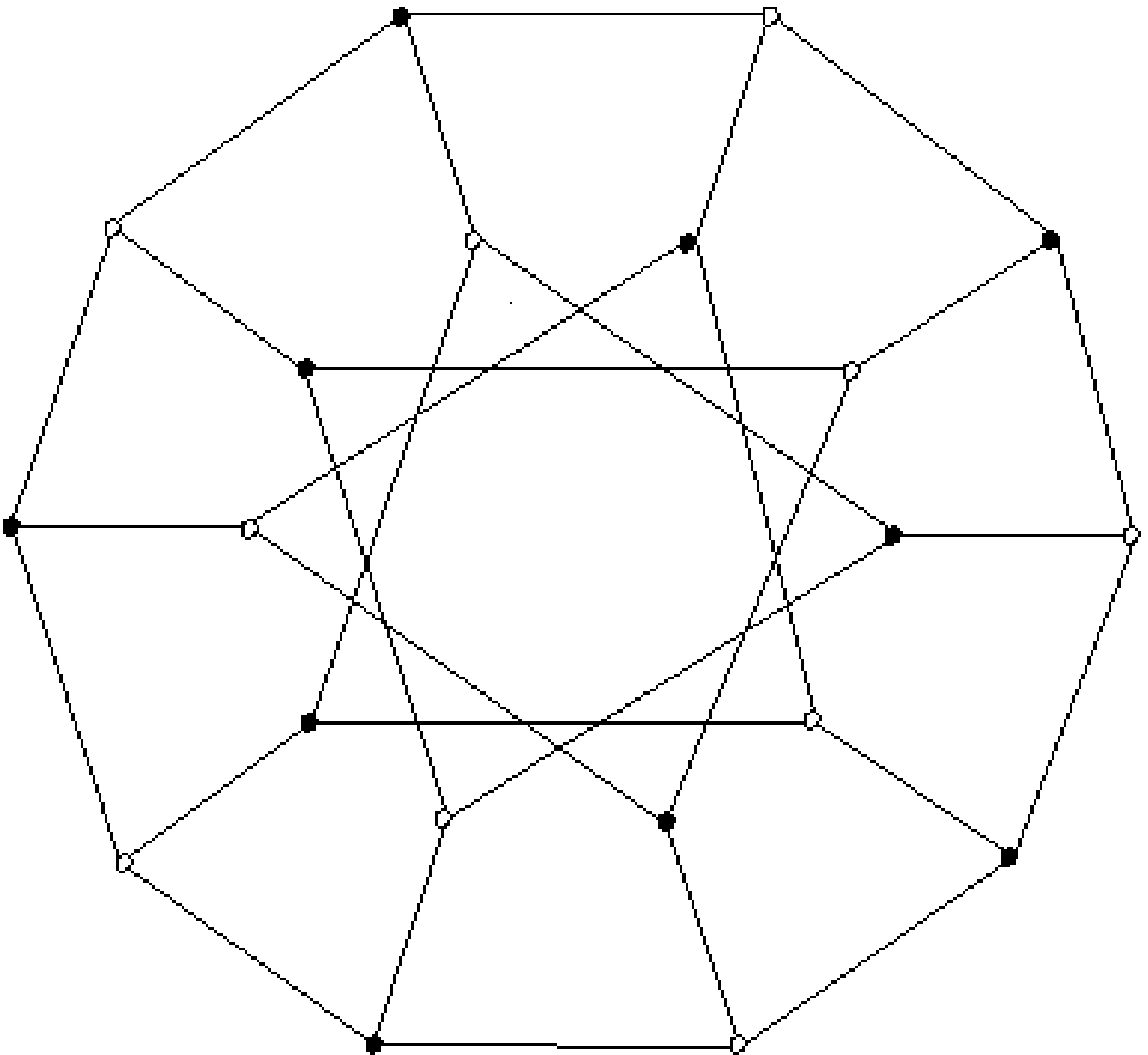}
\caption{ }
\label{desarg}
\end{figure}

\subsection{} All projective planes are symmetric designs with $\lambda = 1$.  Planes for which the Desargues  axiom fails are called {\it non-desarguesian planes}. Some of them can be constructed as projective planes over some non-associative analogs of a finite field  (\emph{near-fields}). In all known examples of projective planes the number of points on a line minus 1 (called the \emph{order} of a plane)  is   always equal to a power of a prime number. The order of the plane $\bbP^2(\bbF_q)$ is of course equal to $q$.

 A non-desarguesian plane of smallest possible order equal to $2^3$ was constructed by Donald Knut (see \cite{Hall}, 12.4).

\section{Configurations in Algebraic Geometry}
\subsection{Linear configurations} A linear $(v_k,b_r)$-configuration over an infinite  field $K$ is a configuration in which the set $\calA$ (resp. $\calB$) is realized by a set of linear subspaces of dimension $d_1$ (resp. $d_2$) of a projective space $\bbP_K^n$ , and the  relation $R$ is an incidence relation  $x\in R(y)$ if $x$ intersects $y$ along a subspace of fixed dimension $s$. The data 
$(\bbP_K^n,s;d_\calA,d_\calB)$ is the type of a linear realization.

For example, any finite linear configuration $\PG(n,r,s;q)$ is a linear configuration over $\bar{\bbF}_q$. One can show that a finite projective plane cannot be realized as a line-point configuration over a field of characteristic zero. I don't know whether it is true for larger $n$.

\smallskip
Note that the sum of two linear configurations is a linear configuration. We consider the corresponding  projective spaces as disjoint subspaces of a projective space and take the joins of subspaces realizing the elements of $\calB$-sets.

\subsection{Dimension reduction}\label{reduction} One can perform the following operations which replace a linear $(v_k,b_r)$-configuration of type $(\bbP_K^n,s;d_\calA,d_\calB)$ with isomorphic  linear $(v_k,b_r)$-configuration. 
\begin{enumerate}
\item [(i)] (Duality): 
The $\bbP_K^n$ is replaced with the dual projective space $\check{\bbP}_K^n$ and the linear subspaces $V$ with the corresponding orthogonal subspace  $V^\perp$.  The type changes to  $(\bbP_K^n,n+s-d_\calA-d_\calB-1;n-d_\calB-1,n-d_\calA-1)$.

\item [(ii)] (Intersection) Take a linear subspace $H$ of dimension $r \ge n-s$ which intersects each subspace from $\calA\cup \calB$ and its intersection subspaces transversally. Here we use the assumption that the ground field is infinite. Then $H$ cuts out a configuration 
of type $(\bbP_K^r,s-n+r, d_\calA-n+r,d_\calB-n+r)$.

\item [(iii)] (Projection) Let $d = \max\{d_\calA,d_\calB\}$. Take a linear subspace $H$ of codimension $c>d+1$ such that for any subspace $V\in \calA$ and any subspace $W\in \calB$ the intersection of the joins
$\langle V,H\rangle$ and $\langle W,H\rangle$ coincide with the join $\langle V\cap W,H\rangle$. Then the projection to $\bbP^{c-1}$ from $H$ defines a configuration of type 
$(\bbP_K^{c-1},s,d_\calA,d_\calB)$. 
\end{enumerate}

 Combining these operations we may sometimes reduce the dimension $n$ of a configuration of type $(\bbP^n,s,d_\calA,d_\calB)$.   Also observe that any abstract $(v_k,b_r)$-configuration  admits a linear realization of type 
$(\bbP_K^r,0,0,r-1)$. It suffices to choose $v$ points in $\bbP^r$ such that any subset of $\calA$ of cardinality $r$  spans a hyperplane. Then for any $B\in\calB$ we consider a hyperplane spanned by $R(B)$ points from $\calA$. This will give a point-hyperplane realization of the abstract configuration.

\subsection{Definition} Let $X$ be a  projective algebraic variety over an algebraically closed field $K$. We say that $X$ realizes a $(v_k,b_r)$-configuration if $X$ contains a set $\calA$ of  $v$ irreducible subvarieties  of dimension $d_\calA$ and a set $\calB$ of $b$ irreducible subvarieties of dimension $d_\calB$ such that   each subvariety from the set $\calA$ intersects $k$ subvarieties from the set $\calB$ along an irreducible   subvariety of dimension $s$ and
 each subvariety from the set $\calB$ intersects $r$ subvarieties from the set $\calA$ along an irreducible subvariety of dimension $s$.

 The type of the realization is the data $(X,s;d_\calA,d_\calB)$.

\subsection{}  It seems that all ``interesting'' classical configurations encountered in algebraic geometry always arise from a  linear configuration when one considers some birational morphism  $X\to \bbP_K^n$.

\subsection{Example}\label{ex1} Any finite linear configuration $\PG(n,r,s;q)$ is realized in  $\bbP_{K}^n$ or in its appropriate blow-up, where $K= \bar{\bbF}_q$.  For example, any projective plane $\PG(2,q)$ is realized in the blowing up of the set of its points. Another linear realization of this configuration can be obtained as follows. The linear system of curves in $\bbP_K^2$ of sufficiently  high degree $d$ passing simply through the set of points $\bbP(\bbF_q)$ will map the blow-up to projective space such that the images of the exceptional curves will be lines and the image of the proper transform of a line in $\PG(2,q)$ will span a linear subspace of codimension $q+1$.

\subsection{Definition} Let $H$ be a subgroup of
$\Sym(\calC)$, where $\calC$ is a $(v_k,b_r)$-configuration. Suppose that $\calC$ is realized in $X$. We say that $H$ is a {\it group of symmetries} of a realization if each element of $H$ is induced  by an automorphism of $X$. The symmetry group of the realization is the largest subgroup of $\Sym(\calC)$  realized as a subgroup of $\Aut(X)$. 

We are interested in  realizations of abstract configurations with as large symmetry group as possible.

\subsection{Example}\label{ex2} In Example \ref{ex1}, any projective automorphism $g\in \PG(n,q)$ is a symmetry of the configuration. It is also a symmetry of the realization. If $k= n-k-1$, then  taking the projective duality, we obtain a switch of the configuration. However, it is not realized as an automorphism of the variety.

\section{Modular configurations}\label{modular} 
\subsection{Abstract $(p(p+1)_p, p^2_{p+1})$-configuration} Let  $p > 2$ be a prime number and $G = (\bbZ/p\bbZ)^2$. We take for $\calA$ the set of cosets of subgroups of $G$ of order $p$ and take for $\calB$ the set $G$. The relation is defined by the inclusion of an element in a coset. There are $p+1$ subgroups of order $p$, each generated either by $(1,0)$ or by $(n,1), n\in \bbZ/p\bbZ$. This defines a $(p(p+1)_p, p^2_{p+1})$-configuration.

More generally, we can consider the set of cosets of subgroups of order $N$ in $G = (\bbZ/N\bbZ)^2$. We have 
\begin{equation}\label{subgroups}
s(N) = N\prod_{p|N}(1+p^{-1}), \quad p\  \textrm{is\  prime}
\end{equation}
Similar to the above we will define an abstract $(Ns(N)_N,N^2_{s(N)})$-configu\-ration.

\subsection{A linear realization}\label{segre} (see \cite{Segre}, \cite{Hulek}) Let $K$ be a field containing $p$ different $p$-th roots of unity. Consider a projective representation of $G$ in $\bbP_K^{p-1}$  (the {\it Schr\"odinger representation}) defined by mapping  $\sigma = (1,0)$ to the projective transformation 
$$\sigma:(x_0,\ldots,x_{p-1})\mapsto (x_{p-1},x_0, x_1,x_2,\ldots,x_{p-2})$$ and mapping 
$ (0,1)$ to the projective transformation 
$$\tau:(x_0,\ldots,x_{p-1})\mapsto (x_0,\epsilon x_1,\ldots,\epsilon^{p-1}x_{p-1}),$$
where $\epsilon$ is a primitive $p$-th root of unity. One checks that each subgroup $H$ of $G$ leaves invariant exactly $p$ hyperplanes. Thus we have a set $\calA$ of $p(p+1)$ hyperplanes. Consider the transformation $\iota$ of $\bbP_K^{p-1}$ defined by the formula
$$\iota:(x_0,x_1,x_2,\ldots,x_{p-2},x_{p-1})\mapsto (x_0,x_{p-1},x_{p-2},\ldots,x_2,x_{1}).$$
The set of fixed points of $\iota$ is equal to the union of a subspace $F^+$ of dimension $\frac{1}{2}(p-1)$ and a subspace $F^-$ of dimension $\frac{1}{2}(p-3)$. The translates $E_g = g(F^-)$ of $F^-$ under the transformations $g\in G$ define a set $\calB$ of $p^2$ subspaces of dimension $\frac{1}{2}(p-3)$. One checks that each subspace $E_g$ is contained in $p+1$ hyperplanes from $\calA$ corresponding to the cosets containing $g$. Also each hyperplane from $\calA$ contains $p$ subspaces from $\calB$ corresponding to the elements in the coset. This defines a linear realization of the $(p(p+1)_p, p^2_{p+1})$-configuration.

\subsection{Hyperflexes of an elliptic curve}
 Let $E$ be an elliptic curve embedded in $\bbP_K^{p-1}$ by the linear system $|px_0|$, where $x_0$ is the zero point.  Each $p$-torsion point $x$ is mapped to a hyperosculating point, i.e. a point $x$ such that there exists a hyperplane which cuts out the divisor $px$ in $E$. The image of the zero point $x_0$ belongs to $F^-$, so that the image of a $g$-translate of $x_0$ is mapped to $E_g$. Thus each hyperosculating point belongs to a unique subspace $E_g$. A hyperplane from $\calA$ cuts out $E$ in $p$ hyperosculating points. When $p =3$ we obtain the famous {\it Hesse $(9_4,12_3)$ configuration} (or {\it Wendepunkts-configuration}) of 9 inflections points of a nonsingular plane cubic which lie by three in 12 lines. 

\subsection{Symmetries} The group $\Gamma = \SL(2,\bbF_p)$ acts naturally on $G$ permuting transitively subgroups of order $p$. There is a projective representation of $\Gamma$ in $\bbP_K^{p-1}$ which leaves the spaces $F^+$ and $F^-$ invariant. The restriction of the representation to each space is derived from a linear irreducible representation of $\Gamma$. The group of symmetries of the modular configuration is isomorphic to the semi-direct product of groups $\SL(2,\bbF_p)$ and $G$. Its order is 
$p^3(p^2-1)$. It is realized by projective transformations of $\bbP_K^{p-1}$. When $p =3$, we get a group of order 216 of projective transformations of the projective plane. It is called the {\it Hesse group}.

\subsection{Modular surfaces}  A modular surface  of level $N > 2$ is an elliptic surface $\pi:S(N)\to B$ representing the universal elliptic curve of level $N$. Its base $B$  is the modular curve  $X(N) = \overline{H/\Gamma(N)}$, where
$$\Gamma(N) = \{A=\begin{pmatrix}a&b\\
c&d\end{pmatrix}\in \SL(2,\bbZ): A\equiv I_2 \quad \textrm{mod}\ N\}.$$
It is a normal subgroup of $\SL(2,\bbZ)$ with quotient group $\SL(2,\bbZ/N\bbZ)$ of order $2\mu(N)$, where
$$\mu(N) = \frac{N^3}{2}\prod_{p|N}(1-p^{-2}), \quad p\  \textrm{is prime }.$$
The Mordell-Weil group $\MW$ of sections of $\pi$ is isomorphic to $(\bbZ/N\bbZ)^2$.  Singular fibres of $\pi$ lie over 
$$c(N) = \mu(N)/N$$ 
cusp points in $X(N)$. Each fibre is of type $\tilde{A}_{N-1}$ (or Kodaira's type $I_N$). For every $b\in B$ let $\Theta_b^{0}$ be the irreducible component of the fibre $F_b = \pi^{-1}(b)$ which intersects the zero section. Then $\Theta_b^{0}$ is intersected by $N$ sections forming a subgroup $H_b$ of order $N$ of  $\MW$. Other components correspond to cosets with respect to $H_b$. We denote by $\Theta_b(\bar{s})$ the component of $F_b$ which is intersected by the sections from a coset $\bar{s} = s+H_b$. For subgroup $H$ let $B_H$ denote the subset of $B$ such that $H = H_b$. The cardinality of this set does not  depend on $H$ and is denoted by $h(N)$.
We have 
\begin{equation}
h(N) = c(N)/s(N) = \frac{N}{2}\prod_{p|N}(1-p^{-1}) = \frac{1}{2}\varphi(N).
\end{equation}
For each subgroup $H$ choose a cusp $b$ such that $H_{b} = H$. Let $\calC'$ be the subset of cusps  obtained in this way. The group $G = \SL(2,\bbZ/N\bbZ)$ acts transitively on the set $\calC$ of cusps. The stabilizer subgroup $G_b$ of a cusp $b$ is a cyclic group of order $N$. Its normalizer $N(G_b)$ is a group of order $Nh(N)$ isomorphic to the stabilizer of the subgroup $H_b$ in the natural action of $\SL(2,\bbZ/N\bbZ)$ on the set of subgroups of $(\bbZ/N\bbZ)^2$.

We have a realization of the $(Ns(N)_N,N^2_{s(N)})$-configuration by a set of irreducible components of fibres $F_b, b\in \calC',$ and the set of sections. The type of this realization is $(S(N),0,1,1)$. We call it a {\it modular configuration} of level $N$. 

\subsection{Hesse pencil} In the case $N = 3$, the modular surface $S(3)$ is obtained by blowing up the base points of  the {\it Hesse pencil} of plane cubic curves
$$\lambda(x_0^3+x_1^3+x_2^3)+\mu x_0x_1x_2 = 0.$$
For each nonsingular member of the pencil the set of its inflection points coincide with the set of base points. The  sections are the exceptional curves of blowing-up. Thus we see that the modular configuration in this case is isomorphic to  the  Hesse $(12_3,9_4)$-configuration.  
 
\subsection{A projective embedding}   Note that the blowing-down morphism $f:S(3)\to \bbP^2$ is given by the linear system $|D|$, where $3D \sim -K_{S(3)}+\sum_{s\in \MW} s $. The canonical class $K_{S(3)} $ of $S(3)$ is equal to $-F$, where $F$ is any fibre of $f$. 
Assume that $N = p > 3$ is prime. One looks for a linear system $|D|$ on $S(p)$ which defines a  birational morphism $f:S(p)\to \bbP^{p-1}$ satisfying the following properties:

\begin{itemize}
\item [(i)] the restriction of $f$ to each  nonsingular fibre $F_b$ is given by the linear system $\bigl|px_b\bigr|$, where $x_b$ is the intersection point of $F_b$ and the zero section $s_0$;
 \item[(ii)] the restriction of $f$ to  each section (identified with the base) is given by the linear system $\bigl|\lambda^{p-3/2}\bigr|$, where $\lambda$ is a divisor of degree $\frac{p^2-1}{24}$ generating the group of $\SL(2,\bbF_p)$-invariant divisor classes on $X(p)$ (see \cite{DO});
\item[(iii)] the map is equivariant with respect to the action of $(\bbZ/p\bbZ)^2\rtimes \SL(2,\bbF_p)$ on $S(p)$ and in $\bbP^{p-1}$.
\end{itemize}
It is known (see \cite{Shioda}) that 
$$K_{S(p)} = 3(p-4)f^*(\lambda),$$ 
$$s^2 =  -\chi(S(p),\calO_{S(p)}) = -p\deg(\lambda)= -\mu(p)/12 = -pc(p)/12.$$
It follows from (i)-(iii) that 
$$pD = \sum_{s\in \MW}s+\frac{1}{2}p(p-1)f^*(\lambda).$$
It is known (see \cite{BH}) that there exists a unique divisor $I$ on $S(p)$ such that $pI\sim \sum s$. Thus we must have $D= I+\frac{p-1}{2}f^*(\lambda)$.
Suppose the linear system $|D|$ contains an invariant subsystem  which defines a birational map $f:S(p)\to \bbP^{p-1}$. Then the image of each section $s$ spans a projective subspace $E_s$ of dimension $(p-3)/2$. The corresponding embedding of $X(p)$  in $E_s$ is the Klein $z$-model of the modular curve (see \cite{DO}). The image of a  nonsingular fibre is a Schr\"odinger-invariant elliptic curve of degree $p$ in $\bbP^{p-1}$. The image of each irreducible component $\Theta_b(\bar{s})$ of a singular fibre is a line. We would like also that the orbit of such a line with respect to the group $N(\SL(2,\bbF_p)_b)$ spans a hyperplane $H(\bar{s})$. Then the configuration formed by the subspaces $E_s$ and hyperplanes $H(\bar{s})$ will be isomorphic to the Segre configuration described in \ref{segre}.   

Let $\Theta_{b_i}(\bar{s}), i = 1,\ldots, (p-1)/2,$ be the components of reducible fibres such that $H_{b_i}$ is the same subgroup $H$ of $\MW$ and $\bar{s}$ is the same coset with of $H$. Then we have to show that 
$$D\sim \sum_{s\in \{s\}}s+\sum_{i=1}^{(p-1)/2}a_i \Theta_{b_i}(\bar{s}),$$
for some positive coefficients $a_i$. Unfortunately, I do not know how to do it for any $p$. For $p = 5$ it is proven in \cite{BH}, and it is also known to be true for $p = 7$.

\section{The Ceva configuration}
\subsection{}  It is is constructed as follows (see \cite{Hirz}). Let $K$ be a field of characteristic prime to $n$ containing the group $\mu_n$ of $n$th roots of unity.  Consider the following $n^3$ points in the projective plane $\bbP_K^2$:
$$P_{0,\alpha} = (0,1,\alpha),\quad  \quad P_{1,\alpha} = (\alpha,0,1),\quad P_{2,\alpha} = (1,\alpha,0),$$
where $\alpha\in \mu_n$. Let $L_{0,\alpha}$ be the line joining the point $(1,0,0)$ with $P_{0,\alpha}$, $L_{1,\alpha}$ be the line joining the point $(0,1,0)$ with $P_{1,\alpha}$, and $L_{2,\alpha}$ be the line joining the point $(0,0,1)$ with $P_{2,\alpha}$. One immediately checks that the lines $L_{0,\alpha}, L_{1,\beta}, L_{2,\gamma}$ intersect at one point $p_{\alpha,\beta,\gamma}$ if and only if $\alpha\beta\gamma = -1$. Thus we obtain $n^2$ points which together with $3n$  lines $L_{i,\alpha}$ form a $(n^2_3,3n_n)$-configuration. It is called the {\it Ceva configuration} and will be denoted by $\Ceva(n)$.

\subsection{Symmetries}For $n \ne 3$, the full symmetry group of this configuration is equal to the semi-direct product $\mu_n^2\rtimes S_3$. It is realized by the group of projective transformations of $\bbP_K^2$ generated by permutation of the coordinates and homotethies $(x_0,x_1,x_2)\mapsto (\alpha x_0,\beta x_1,\gamma x_2)$, where $\alpha\beta\gamma = 1$. 

When $n = 3$, the symmetry group is larger. If one realizes the configuration in $\bbP^2(\bbF_4)$, then we get an additional symmetry realized by the Frobenius automorphism. Also we have a duality automorphism. The resulting group is the Hesse group of order 216. 

\subsection{Blow-up}Blowing up the set of points $p_{\alpha,\beta,\gamma}$, we get a rational surface $V$ together with a morphism $\pi:V\to \bbP^1$ whose general fibre is a nonsingular curve of genus $g = (n-1)(n-2)/2$. There are 3 singular fibres; each is the union of $n$ smooth rational curves   with self-intersection $1-n$ intersecting at one point. The morphism admits $n^2$ disjoint sections; each is a smooth rational curve with self-intersection $-1$. The Ceva configuration is realized by the set of sections and the set of irreducible components of singular fibres. If $n\ne 3$, the symmetry group of the configuration is realized by an automorphism group of the surface. As we  will see  in the last section there is a realization of $\Ceva(3)$ which realizes a subgroup of index 2 of $\Sym(\Ceva(3))$. 

In the case $n = 2$, $V$ is a minimal nonsingular model of a 4-nodal cubic surface, and the morphism is induced by a pencil of conics. In the case $n = 3$, $V$ is a rational elliptic surface with 3 singular fibres of type $IV$ with Mordell-Weil group containing a subgroup isomorphic to $(\bbZ/3\bbZ)^2$.

\section{$v_3$-configurations}

\subsection{Linear realization} All $v_3$-configurations can be linearly realized in $\bbP_K^3$ by $v$ points in general linear position and $v$ planes.

\subsection{Symmetries} The group of symmetries of an abstract $v_3$-configu\-ration with $s$-regular  Levi graph is of order $2^s3n$ (see \cite{Cox}).

\subsection{$v\le 6$} There is a unique $v_3$-configuration for $v = 4,5,6$. 

A $4_3$-configuration is realized by the vertices and faces of a tetrahedron. Its Levi graph is 1-regular 

A $5_3$-configuration can be realized by 5 points $p_1,\ldots,p_5$ in general linear position in $\bbP_K^3$ and the planes 
$$\langle p_1,p_2,p_5\rangle,\   \langle p_1,p_2,p_3\rangle,\  \langle p_2,p_3,p_4\rangle, \ \langle p_3,p_4,p_5\rangle,\  \langle p_1,p_4,p_5\rangle.$$
Its Levi graph is 2-regular.

A $6_3$-configuration can be realized by 6 points $p_1,\ldots,p_6$ in general linear position in $\bbP_K^3$ and the planes 
$$\langle p_1,p_3,p_4\rangle,\   \langle p_1,p_4,p_6\rangle, \ \langle p_1,p_2,p_5\rangle, \ \langle p_2,p_3,p_6\rangle, \ \langle p_2,p_4,p_5\rangle,\  \langle p_3,p_5,p_6\rangle$$
Its Levi graph is 2-regular.

All these configurations satisfy the property that there exist two distinct  $x,y\in\calA$ such that $\# R(x)\cap R(y) > 1$. Obviously such configuration cannot be realized by points and lines in projective space. 

 \subsection{$7_3$-configuration} There is a unique  line-point $7_3$-configuration. It is isomorphic to the Fano plane. The symmetry group of its natural realization in $\bbP_{\bar{\bbF}_2}$ is isomorphic to $\PGL(2,\bbF_2)$ and is smaller than the group of symmetry of the configuration. The latter is generated by this subgroup, the Frobenius and  a switch.  Let $X$ be the blow-up of this realization. Now the set $\calA$ is represented by seven $(-1)$ curves and the elements of $\calB$ are realized by seven $(-2)$-curves. The linear system of cubics through the seven points is generated by the cubics 
$$x_0x_1(x_0+x_1) = 0, \quad x_0x_2(x_0+x_2) = 0,\quad x_1x_2(x_1+x_2) = 0.$$
It defines a regular map of $X$ to itself by interchanging the sets $\calA$ and $\calB$. Note that the map is of degree 2 and inseparable. So it is `almost a switch' since it is bijective on the set of points $X(\bbF_2)$. We shall discuss later a different realization of the Fano plane by curves on a K3 surface which realizes a switch.

It is known the Fano plane  can be  realized by lines and hyperplanes in $\bbP^4(\bbR)$ (\cite{Kelly}). It is not known whether this is possible for other finite projective 
planes.

\subsection{$8_3$-configuration}  There is a unique abstract  line-point $8_3$-configu\\
ration. It is called the {\it M\"obius-Kantor} configuration. Here is its Levi graph (see \cite{Cox}):

\begin{figure}[h]
\begin{center}
\includegraphics[width=2in]{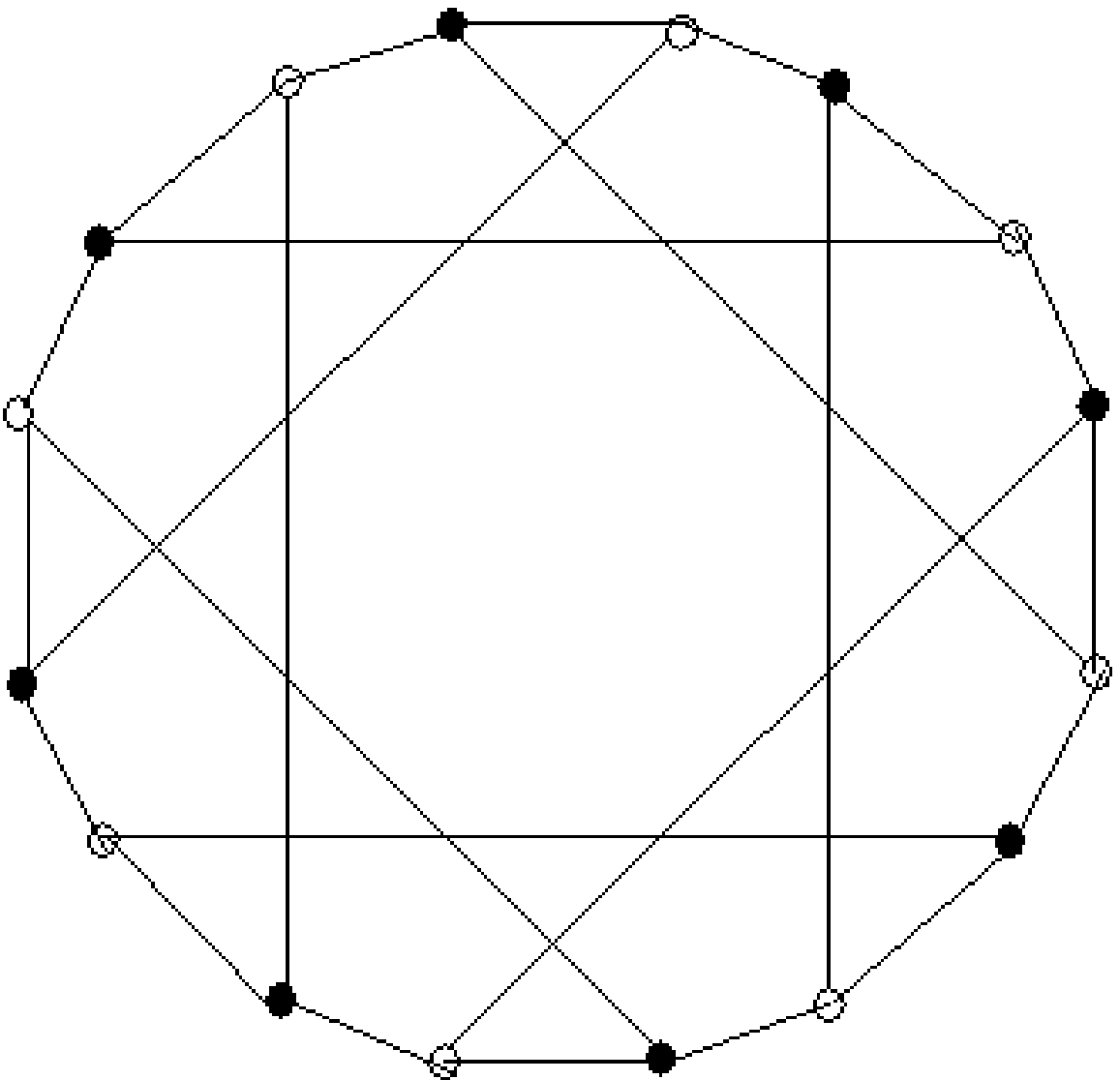}
\caption{ }
\label{moebius}
\end{center}
\end{figure}
A linear realization of the M\"obius-Kantor configuration can be obtained from the Hesse $(12_3,9_4)$-configuration by deleting one inflection point and the four lines to which it is incident (see \cite{Levi}).

\subsection{$9_3$-configurations} There are three non-isomorphic  line-point abstract $(9,3)$-configurations (see \cite{Hilbert},\cite{Levi}). All of them are realized by points and lines in the plane.

 The first $9_3$-configuration is the {\it Brianchon-Pascal configuration}. 
\begin{figure}[h]
\begin{center}
\includegraphics[width=2.4in]{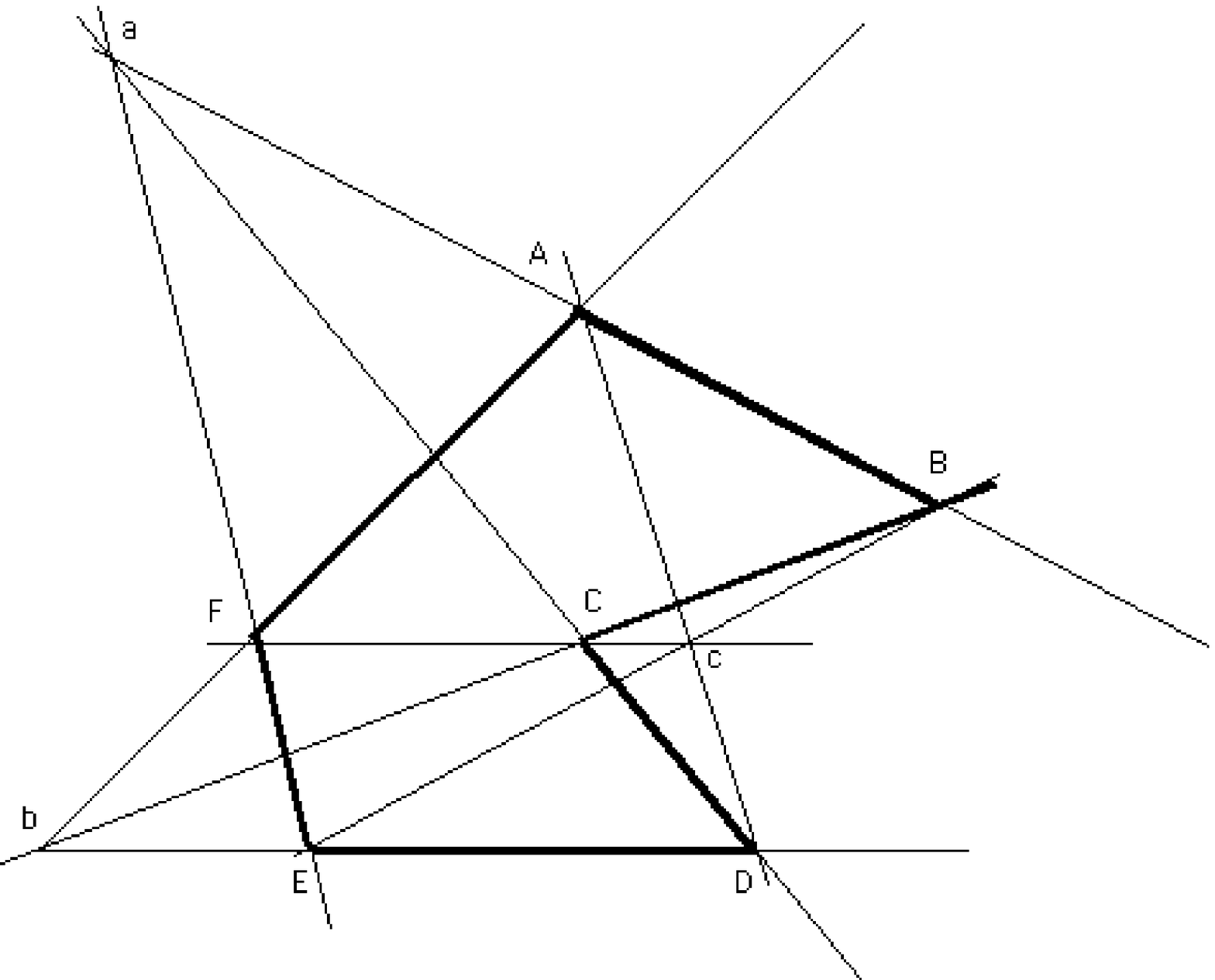}
\caption{ }
\label{brianchon}
\end{center}
\end{figure}

It admits a realization  by three diagonals and six sides of a  hexagon $ABCDEF$ such that the sides $AB,CD,EF$ intersect at a point $v$ and the sides $BC,DE,FA$  intersect at a point $b$.  The Brianchon Theorem asserts that the diagonals intersect at a point $c$. This gives us 9 lines and nine points (six vertices and  the points $a,b,c$).

The Brianchon-Pascal configuration is isomorphic to its dual  realization:

\begin{figure}[h]
\begin{center}
\includegraphics[width= 2in]{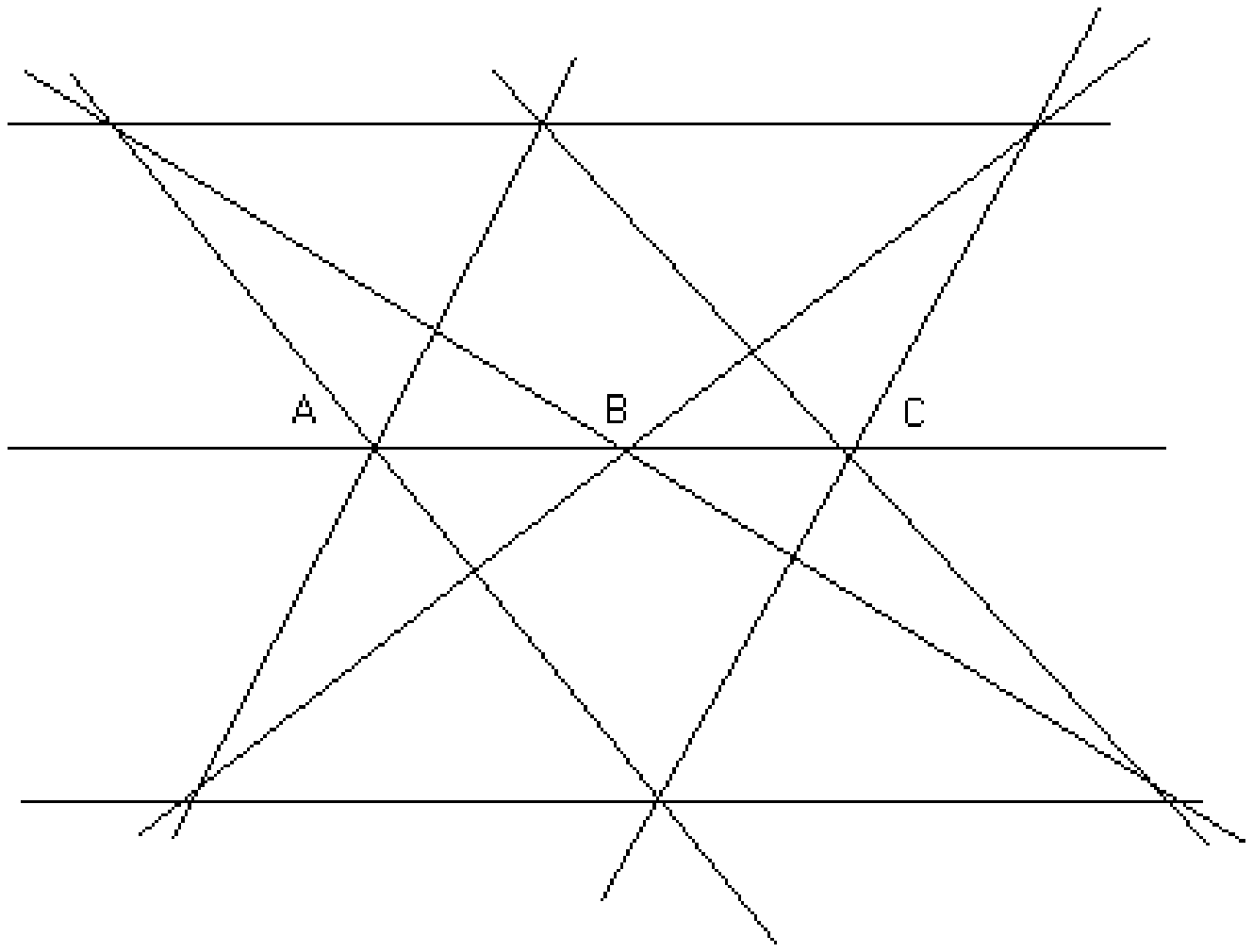}
\caption{ }
\label{ }
\end{center}
\end{figure}
The {\it Pascal  Theorem} asserts that the points $A,B,C$ are on a line (the {\it Pascal 
line}).

The Brianchon configuration is isomorphic to the Ceva configuration for $n = 3$. So, its symmetry group is isomorphic to the Hesse group.

Here are the pictures of other two $9_3$-configurations.
\begin{figure}[h]
\begin{center}
\includegraphics[width=2.5in]{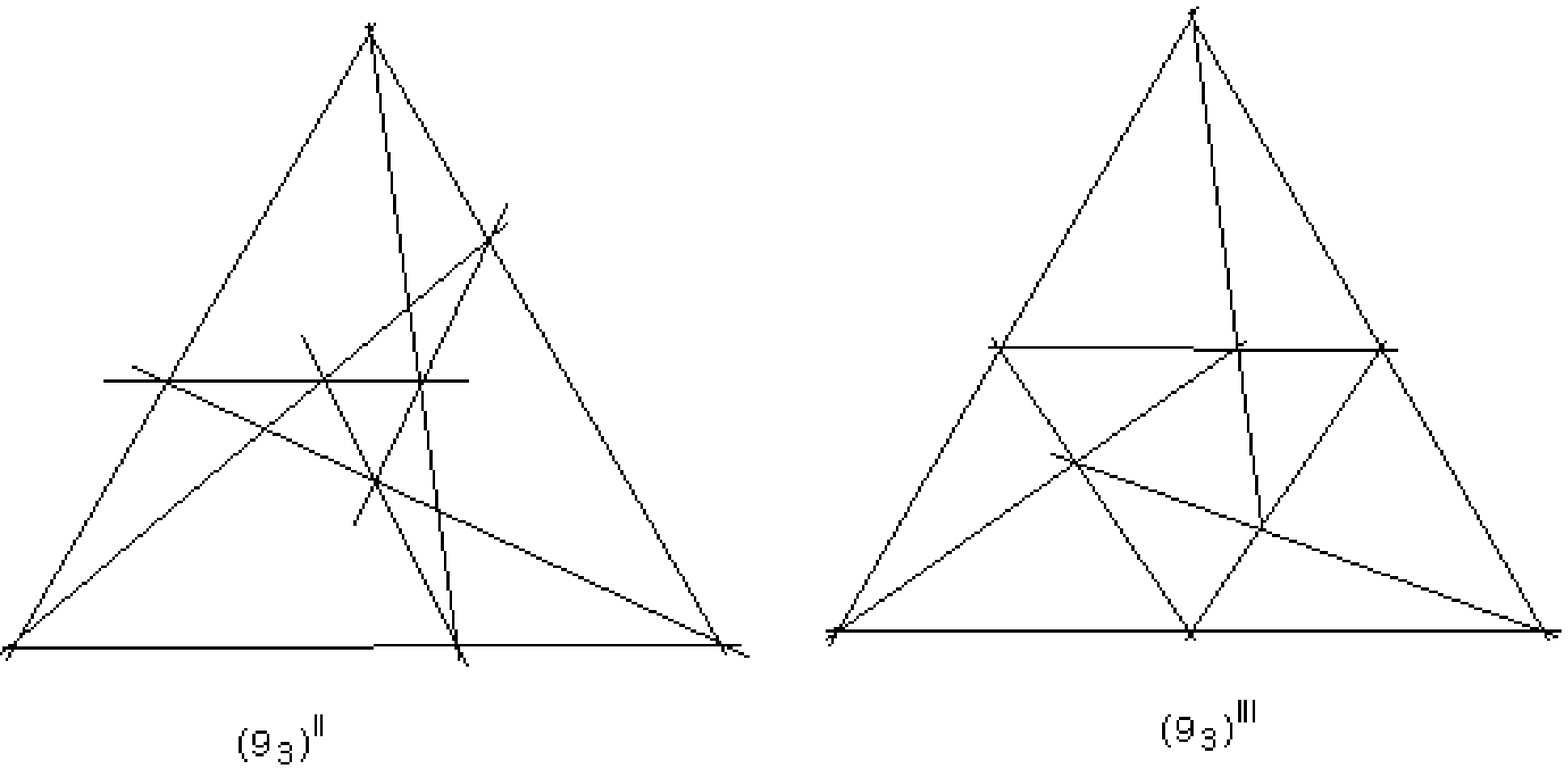}
\caption{ }
\label{ }
\end{center}
\end{figure}
\newpage

To see that these three configurations are not isomorphic one argues as follows. For any configuration let $S$ be the set of pairs $x,y\in \calA$ such that $R(x)\cap R(y) = \emptyset$. We say that two elements $x,y$ in $\calA$ are $S$-equivalent  if either $x = y$ or there exists a  sequence of pairs $s_1 = (x,x_1), s_2 = (x_1,x_2),\ldots, s_n = (x_{n-1},y)$ in $S$.   Clearly this is an equivalence relation. Also it is clear that any isomorphism of configurations sends an $S$-equivalence class to an $S$-equivalence class.

In the examples above, the first configuration has three $S$-equivalence classes, each consists of 3 elements. The second configuration consists of one $S$-equivalence class. The third  configuration has 2 equivalence classes, of 6 and 3 elements. 

The first two configurations are regular. The configuration of third type is obviously non-regular (in a regular configuration all $S$-equivalence classes consist of the same number of elements).

\subsection{$10_3$-configurations} We have already described a realiziation of  a Desargues $10_3$-configuration  by points and lines on a projective plane over a field of arbitrary characteristic.   Another realization of this configuration is as follows. It is known that a general cubic surface can be written as a sum of 5 cubes of linear forms. The locus of the linear forms form the Sylvester pentahedron associated to the cubic surface. The Hessian  surface is a quartic which contains the 10 edges of this pentahedron (i.e. the lines of intersection of pairs of the planes) and has 10 nodes corresponding to the vertices of the pentahedron. A minimal resolution of the Hessian quartic is a K3 surface with two sets of ten $(-2)$-curves forming  a $10_3$-configuration isomorphic to the Desargues configuration.

The Brianchon  Theorem admits another version which  asserts that the diagonals of a hexagon circumscribed about an irreducible  conic intersect at one point. Its dual Pascal theorem asserts that three pairs of opposite sides of a hexagon inscribed in an irreducible conic intersect at three collinear points. 

Both versions of the Brianchon Theorem can be seen by projecting a hexagon formed by rulings of a quadric from either a point on the quadric or from a point outside the quadric (see \cite{Hilbert}).

The corresponding pictures are called  {\it Pascal Hexagrams}:

\begin{figure}[h]
\begin{center}
\includegraphics[width=2.3in]{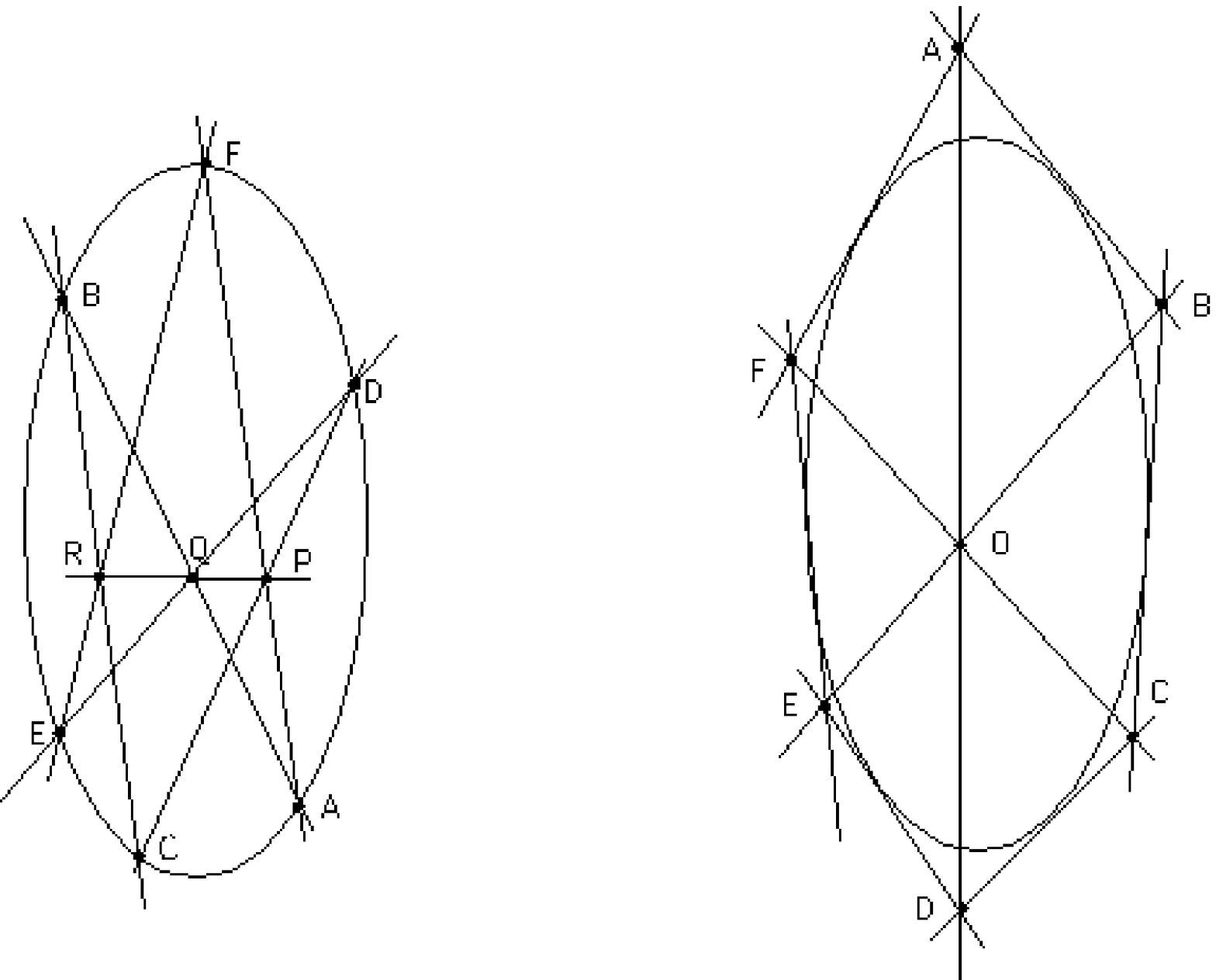}
\caption{ }
\label{hexagram}
\end{center}
\end{figure}

Thus if we chose a hexagon which is inscribed in a conic $C_1$ and  also circumscribed about another conic $C_2$ (the conics are called Poncelet $6$-related), we can apply  both theorems to obtain a $10_3$-configuration isomorphic to the Desargues configuration.

There are 10 non-isomorphic abstract line-point $10_3$-configurations. Only one of them cannot be realized in the complex projective plane. This result is due to S. Kantor \cite{Kantor} (see a modern proof in  \cite{Betten}). Here is a  regular self-dual $10_3$-configuration which is not isomorphic to the Desargues configuration.
\begin{figure}[h]
\begin{center}
\includegraphics[width=1.5in]{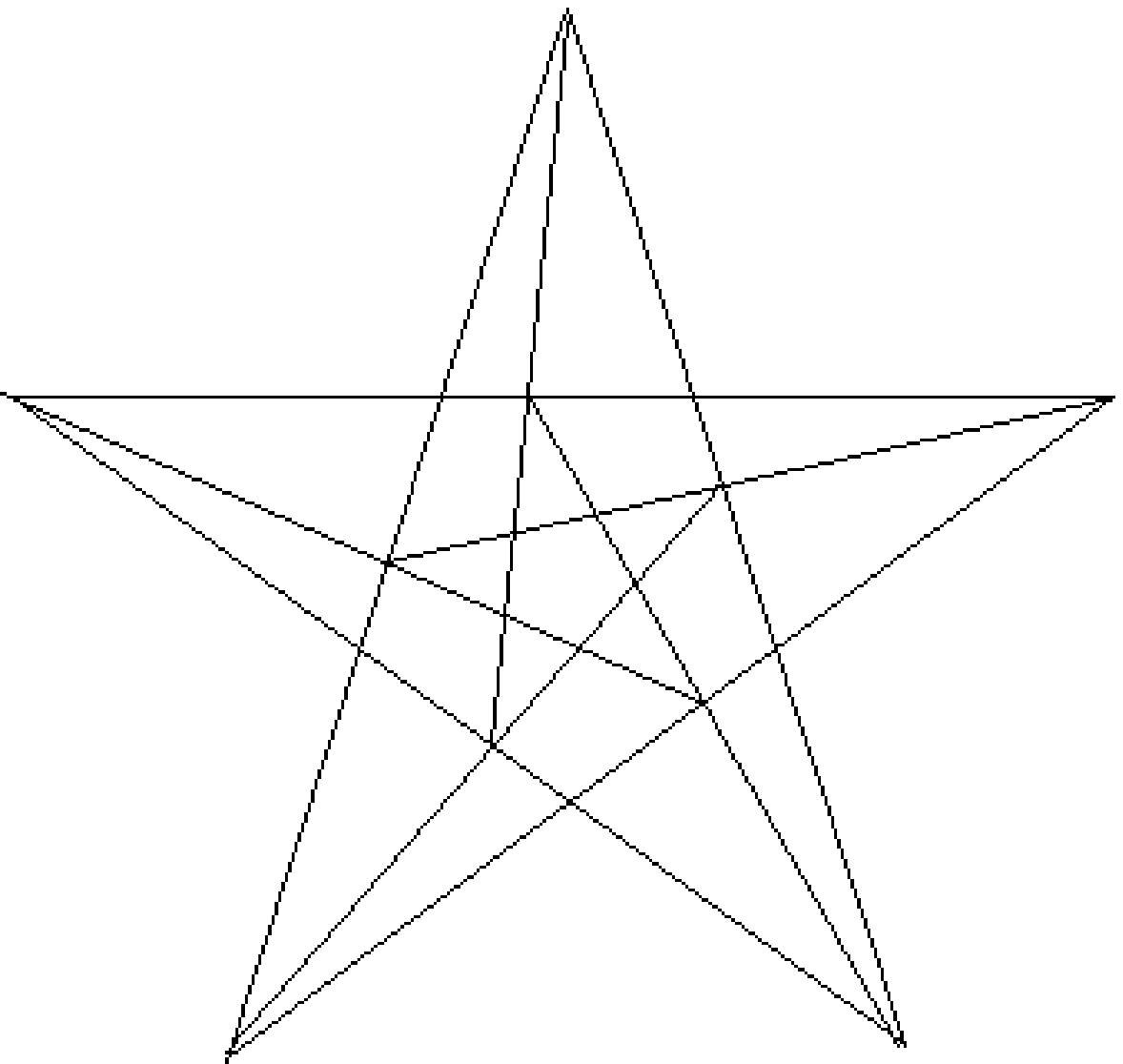}
\caption{ }
\label{ }
\end{center}
\end{figure}

The  number of non-isomorphic abstract line-point  $(11_3)$-configurations  is equal to 31 (\cite{Mart}, \cite{Daub1}) and  for $(12_3)$-configurations is equal to 228 (\cite{Daub2}).

\section{The Reye $(12_4,16_3)$-configuration}
\subsection{The octahedron configuration} The first construction is due to T. Reye \cite{Reye}, p.234. Consider a cube in $\bbR^3$. Let $\calA$ be the set of 16 lines in $\bbP^3(\bbR)$ which consists of 12 edges and  4 diagonals of the cube. Let $\calB$ be the set of 12 points in $\bbP^3(\bbR)$ which consists of 8 vertices of the cube and 4 vertices of the tetrahedron formed by the intersection point of the diagonals and three points at infinity where the edges of the cube intersect. Each point from $\calB$ belongs to four lines from $\calB$. Each line from $\calB$ contains three points from $\calB$. This defines a $(12_4,16_3)$-configuration.

\begin{figure}[h]
\begin{center}
\includegraphics[width=1.5in]{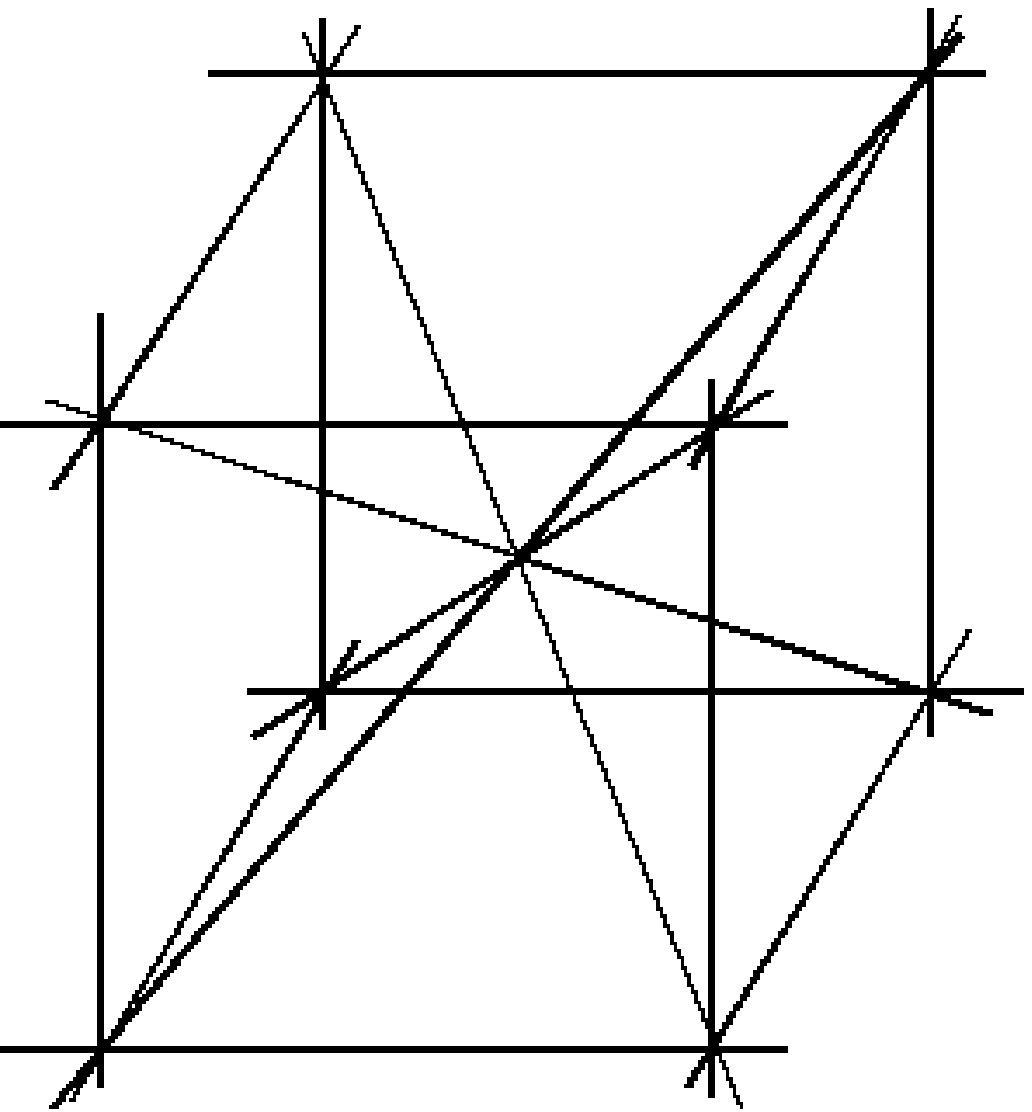}
\caption{ }
\label{ }
\end{center}
\end{figure}
\subsection{A net of diagonal quadrics}  Here is another way to define this configuration. Consider a net $\calN$ of quadrics in $\bbP^3(\bbC)$ spanned by three diagonal quadrics 
\begin{equation}\label{quadrics}
Q_1:\sum_{i=0}^3x_i^2 = 0,\quad Q_2:\sum_{i=0}^3a_ix_i^2,\quad Q_3:\sum_{i=0}^3b_ix_i^2 = 0.
\end{equation}
We assume that all $3\times 3$ minors of the matrix
\[\begin{pmatrix}1&1&1&1\\
a_0&a_1&a_2&a_3\\
b_0&b_1&b_2&b_3\end{pmatrix} \]
are nonzero. 

Let $H: =\det (\lambda Q_1+\mu Q_2+\gamma Q_3) = 0$ be the discriminant curve of the net $\calN$. It consists of four lines in general linear position. The six intersection points of the lines correspond to 6 reducible quadrics in $\calN$. Each such a quadric consists of two planes intersecting along an edge of the coordinate tetrahedron $T$. Let $A_1,\ldots,A_4$ be the vertices of $T$ and let   $Q_{ij} = P_{ij}\cup P_{ij}'$ be the reducible quadric from $\calN$ containing the edge $\langle A_i,A_j\rangle $. The set of 8 base points of the net lie by four in each plane $P_{ij}, P_{ij}'$. Let $\langle A_i,A_j\rangle $ and $\langle A_k,A_l\rangle $ be two opposite edges. The quadric $Q_{ij}$ intersects $\langle A_k,A_l\rangle $ at two points $p_{kl},p_{kl}'$ and the quadric $Q_{kl}$ intersects $\langle A_i,A_j\rangle $ at two points $p_{ij},p_{ij}'$. The eight base points of $\calN$ lie by two on the four lines $\langle p_{ij},p_{kl}\rangle , \langle p_{ij},p_{kl}'\rangle , 
\langle p_{ij}',p_{kl}\rangle , \langle p_{ij}',p_{kl}'\rangle $. The three pairs of opposite edges give us  12 lines which together with eight base points form a $(8_3,12_2)$-configuration isomorphic to the configuration of vertices and edges of a parallelepiped $\Pi$. Let $Q_{ij},Q_{ik},Q_{il}$ be three reducible quadrics from $\calN$ containing a vertex $A_i$. Each pair intersects along the same set of 4 lines which join $A_i$ with two base points. Altogether we find 16 lines which contain by three the 12 singular points of the cubic complex of lines  and form a $(12_4,16_3)$-configuration isomorphic to the octahedron configuration.

\subsection{The Hesse-Salmon configuration} Consider now the projection of the Reye configuration to the plane from a general point in the space. We get a line-point $(12_4,16_3)$-configuration in the plane. We will show that it is isomorphic to the following {\it Hesse-Salmon configuration} (\cite{Hesse},\cite{Salmon}, \cite{Zah},\cite{Cox}). Let $C$ be a nonsingular plane cubic curve. Take a line which intersects $C$ at three points $a,b,c$. The pencil of lines through each point $p_i$ contains 4 tangent lines  to $C$. Let 
$$(A_1,A_2,A_3,A_4), \quad (B_1,B_2,B_3,B_4),\quad (C_1,C_2,C_3,C_4),$$
 be   the three sets of four  points of contact. Fix a group law on $C$ by taking one of the inflection points as the zero point. Then 
\[a+b+c = 0,\ 2A_i+a = 0, \  2B_i+b = 0,\  2C_i+c = 0,\  i = 1,2,3.\]
We can write
\[A_2 = A_1+\alpha,\quad  A_3 = A_1+\beta,\quad   A_4 = A_1+\alpha+\beta,\]
where $\alpha,\beta$ are two distinct 2-torsion points on $C$.
After reordering, if needed, we can  write
\[B_2 = B_1+\beta, \quad B_3 = B_1+\alpha+\beta, \quad B_4 = B_1+\alpha,\]
\[C_2 = C_1+\alpha+\beta, \quad C_3 = C_1+\alpha, \quad C_4 = C_1+\beta.\]
Hence, we obtain that
{\small \begin{alignat}{3}\label{lines}
A_1+B_1+C_1&=  A_1+B_2+C_4 & = A_1+B_3+C_2& = A_1+B_4+C_3 &=    &0\\ \notag
A_2+B_1+C_3& = A_2+B_2+C_2& = A_2+B_3+C_4& = A_2+B_4+C_1 &=   &0,\\\notag
A_3+B_1+C_4 &= A_3+B_2+C_1 &= A_3+B_3+C_3& = A_3+B_4+C_2 &=  & 0,\\\notag
 A_4+B_1+C_2&=   A_4+B_2+C_3 & = A_4+B_3+C_1& = A_4+B_4+C_4 &=   &0\\\notag
\end{alignat}}
This shows that the 12 contact points lie on 16 lines forming a $(12_4,16_3)$-configuration.

\subsection{}Let us now  show that the Reye configuration is isomorphic to the Hesse-Salmon  configuration. We have to recall a few constructions from the theory of linear systems of quadrics. Let $\calN$ be the net of quadrics spanned by the quadrics from \eqref{quadrics}. Let $G_{1,3}$ be the Grassmannian of lines in $\bbP^3$ and 
$X$ be the set of lines  $l\in G_{1,3}$ such that there exists $Q\in \calN$ containing $l$.  Let 
$D_i$ be the divisors of bi-degree  $(1,1)$ in $\bbP^3\times \bbP^3$ defined by the  symmetric bilinear forms associated to the quadrics $Q_i$. The intersection $D_1\cap D_2\cap D_3$ is a 3-fold $Y$ in $\bbP^3\times \bbP^3$. The variety $X$ is the closure of the set of lines $\langle x,y\rangle$, where  $(x,y)\in Y, x\ne y.$ Let $\pi$ be a general plane in $\bbP^3$. Restricting quadrics from $\calN$ to $\pi$ we obtain a net of conics in $\bbP^2$. It is easy to see that set of lines in $\pi$ contained in a conic from this net is a nonsingular plane cubic in the dual plane. Thus the number  of lines in $\pi$ passing through a fixed general point is equal to 3. This shows that $X$ is a cubic complex in $G_{1,3}$, i.e., it is cut out by a cubic hypersurface in the Pl\"ucker embedding of $G_{1,3}$.  For any $x\in X$ let $K(x)$ be the set of lines in $X$ which contain $x$. This is a plane section of $X$, and hence  is a plane cubic curve. Thus the union of lines from $K(x)$ form a cubic cone with vertex at $x$. A point $x\in \bbP^3$ is called a {\it singular point} of a line complex if each line through this point belongs to the complex. In our cubic complex $X$ there are 12 singular points. They are the eight base points of the pencil and the four vertices of the coordinate tetrahedron $T$.

 Now consider projection from a general point $x$ in $\bbP^3$. The  lines in $K(x)$ project to a nonsingular plane cubic $C$.  Since $K(x)$ contains the lines joining $x$ with a singular point of $X$, we see that the projections of the 12 singular points lie on $C$. Let $l_1,l_2,l_3,l_4$ be four lines among the sixteen lines which contain the vertex $A_1$. Let $B_1,C_1$ be the base points lying on $l_1$, $B_2,C_4$ be the base points lying on $l_2$, $B_3,C_2$ be the base points lying on $l_3$, and $B_4,C_3$ be the base points lying on $l_4$. We assume that no two  of the $C_i$'s and  no two of the $B_i$'s lie on the same edge of the parallelepiped $\Pi$. Then we may assume that the lines $\langle B_1,C_3\rangle,\langle B_2,C_2\rangle, \langle B_3,C_4\rangle, \langle B_4,C_1\rangle$ intersect at $A_2$, the lines $\langle B_1,C_4\rangle, \langle B_2,C_1\rangle, \langle B_3,C_3\rangle, \langle B_4,C_2\rangle$ intersect at $A_3$,  and the lines $\langle B_1,C_2\rangle, \langle B_2,C_3\rangle, \langle B_3,C_1\rangle, \langle B_4,C_4\rangle$ intersect at $A_4$. 

If we use the group law on the cubic $C$, we may assume that equations \eqref{lines} are satisfied. 
Now the equations 
$$A_1+B_3+C_2 = A_1+B_2+C_4 = A_2+B_2+C_2 = A_2+B_3+C_4 = 0$$
imply that $2A_1 = 2A_2$, and similarly we get $2A_1 = 2A_3 = 2A_4$. This shows that the tangent lines to $C$ at the points $A_1,A_2,A_3,A_4$ intersect at some point $a\in C$. In the same way we obtain that the tangent lines at the points $B_i$ intersect at some $b\in C$ and the tangent lines at the points $C_i$ intersect at some point $c\in C$. We get $2A_1 + a = 0, 2B_1 +b = 0, 2C_1 + c = 0$, and the first equation in \eqref{lines} gives $a+b+c = 0$, i.e. the points $a,b,c$ are collinear. This proves that the projection of the Reye configuration is equal to the Hesse-Salmon configuration.

\subsection{Symmetries} The group of abstract symmetries is isomorphic to the octahedron group of order $24$.  In is realized in the octahedron realization. The group of symmetries of the net of quadrics realizations is isomorphic to the group of order 8 generated by projective 
transformations $x_i\mapsto \pm x_i$.  The group of symmetries of the Hesse-Salmon realization is trivial.

\section{$v_{v-1}$-configurations} 

\subsection{} A $v_{v-1}$-configuration is obviously unique and is not a line-point abstract configuration. The proper symmetry group of this configuration is isomorphic to the permutation group $S_a$. There is also a duality defined by switching a point $x\in\calA$ with the unique point $y$ from $B$ such that $x\not\in R(y)$. The configuration can be  realized by points and hyperplanes in $\bbP_K^{a-1}$ by taking for $\calA$ a set of $v$ points in general position and taking for $\calB$ the set of hyperplanes spanned by all points except one.

\subsection{The double-six} This is $(6_5)$-configuration realized by a  double-six of lines on a nonsingular cubic surface.  The full group of symmetry of this configuration is the double extension of $S_6$. It is generated by permutation of lines in one family and a switch. In the full group of symmetries of 27 lines on a cubic surface this is the subgroup of $W(E_6)$ of index 36 which fixes a subset $\{\alpha,-\alpha\}$, where $\alpha$ is a positive root. One can realize the full symmetry group over a field of characteristic 2 by considering the Fermat cubic surface
$$x_0^3+x_1^3+x_2^3 +x_3^2 = 0.$$
Its automorphism group is isomorphic to the Weyl group $W(E_6)$ (see \cite{Hirsch}), Theorem 20.3.1).

\subsection{The Fano $(10_9)$-configuration} Let $X$ be an Enriques surface over an algebraically closed field $K$. Assume that $X$ has no smooth rational curves. Then $X$ can be embedded in  $\bbP^5$ as a surface of degree 10 (a {\it Fano model}, see \cite{DoRe}). The surface $X$ contains 10 elliptic pencils $|2F_i|$ such that $F_i\cdot F_j = 1$ if $i\ne j$.  The divisor class of a hyperplane section $H$ of $X$
 satisfies $3H\equiv F_1+\ldots+F_{10}$. 

Assume that the characteristic of $K$ is not equal to $2$. Each elliptic fibration contains two double fibres $2F_i$ and $2F_i'$. The image of the curves $F_i, F_i$'s are plane cubics. Let us choose one of the double fibres in each fibration and denote by $\calA$ the set of the planes in $\bbP^5$ spanned by the corresponding fibres.  We take for $\calB$ the set of planes corresponding to the remaining 10 double fibres. We get a $(10_9)$-configuration formed by the sets $\calA$ and $\calB$ and the incidence relation defined by $(\pi,\pi')\in R$ if $\pi\cap \pi'$ is a one-point set. 
 
If $K$ is of characteristic 2 and $K_X = 0$, then $F_i = F_i'$. In this case we  take $\calA = \calB$.

In general, the automorphism group of the realization is trivial. It is an interesting problem to find $X$ realizing non-trivial symmetries of the configuration.

\subsection{Coble's configurations} A generalization of a double-six is due to A. Coble (\cite{Coble2}). It is a $(v_{v-1}$)-configuration, where  $v = N_n = \binom{n+2}{2}$. One takes a set $\Sigma$ of $N_n$ general points $\Sigma =\{p_1,\ldots,p_{N_n}\}$ in the projective plane over an algebraically closed field of any characteristic. For each points $p_j$ there is a unique smooth curve $C_j$ of degree $n$ which passes through the set $\Sigma\setminus \{p_j\}$.  Let $X_n$ be the blow-up of $\bbP^2$ along the subset $\Sigma$. We take $\calA$ to be the set of exceptional curves $E_j$, and $\calB$ to be the set of proper transforms $\bar{C}_j$ of the curves $C_j$. Note that 
$$\bar{C}_j^2 = \bar{C}_j\cdot K_{X_n} = n(n-3)/2.$$ 
 Note that the linear system of curves of degree $n+1$ through the set $\calA$  maps $X$ isomorphically to a surface $W_n$ of degree $N_{n-1}$ in $\bbP^{n+1}$.  It is called a {\it White surface}. Clearly, $W_1 \cong \bbP^2$ and $W_2$ is a cubic surface. The image of the proper inverse transform of $C_j$ in $X_n$ is mapped to a curve in $\bbP^{n+1}$ lying in a linear subspace of dimension $n-1$. The image of each exceptional curve is a line in $W_n$. Thus our configuration is isomorphic  to a configuration of lines and subspaces of dimension $n-1$ in $\bbP^{n+1}$. 
Unfortunately, when $n > 2$, no symmetry of our configuration is  realized by an automorphism of the surface $W_n$.

\section{The Cremona-Richmond $15_3$-configuration} 
\subsection{} This configuration  can be realized in $\bbP_K^4$ as follows (see \cite{Rich}). Take 6 points $p_1,\ldots,p_6$ in general linear position. For any subset $I$ of the set $[1,6] = \{1,\ldots,6\}$ let $\Pi_I$ be the linear subspace spanned by the points $p_i, i\in I$. We have 15 ``diagonal'' points $p_{I} = \Pi_I\cap \Pi_{\bar {I}}$, where $\#I = 2,$ and the bar denotes the complementary set. The points $p_{I},p_{J},p_{K}$ with $[1,6] = I\cup J\cup K$ lie on the ``transversal line'' $\ell_{I,J,K}$ , the intersection line of the hyperplanes $\Pi_{\bar{I}},\Pi_{\bar{J}},\Pi_{\bar{K}}$. The 15 diagonal points and 15 transversal lines form a $15_3$-configuration. 
By projecting from a general line, we get a line-point plane $15_3$-configuration. 

\subsection{The dual configuration} Let us  identify  the dual projective space $\check{\bbP}^4$ with the hyperplane $V\cong \bbP^4$ defined by the equation $x_0+\cdots+x_5 = 0$ in $\bbP^5$.  Let the dual of points  $p_i$ to be the  hyperplanes in $V$ defined by the additional equations $x_i = 0$.  Then the dual of the diagonal points are the  hyperplanes 
$H_{ab}$ defined by the additional equations $x_a +x_b = 0$, and the dual of transversals  are the planes $\pi_{ab,cd,ef}$ defined by the additional equations $x_a+x_b = x_c+x_d = x_e+x_f = 0$. Let $S_6$ act in $V$ by permuting the coordinates $x_i$. This is the projective representation of $S_6$ defined by the natural 5-dimensional irreducible representation  of $S_6$. The dual of the diagonal points are the hyperplanes of fixed points of 15 transpositions $(ab)$. The dual of the transversals are the planes of fixed points of 15 permutations of cycle type $(ab)(cd)(df)$. The group $S_6$ is isomorphic to the proper symmetry group of the arrangement via its natural action on itself by conjugation. The full symmetry group of the arrangement is isomorphic to the group of outer automorphisms of $S_6$.

\subsection{The Segre cubic} The Segre cubic $\calS_3$ is defined by the equation 
$x_0^3+\cdots+x_5^3 = 0$ in $V$. It admits a nice moduli interpretation (see \cite{Coble}, \cite{DoOr}). It has 10 ordinary double points with coordinates $(1,1,1,-1,-1,1)$ and their permutations. Each plane $\pi_{ab,cd,ef}$ is contained in $\calS_3$. The hyperplane $H_{ab}$ cuts out in $\calS_3$ the union of three planes $\pi_{ab,cd,ef}, \pi_{ab,ce,df},\pi_{ab,cf,de}$. Each plane $\pi_{ab,cd,ef}$ contains 4 singular points, and each singular point is contained in 6 planes. This defines a $(15_4,10_6)$-configuration . 

\subsection{Cubic surfaces}\label{cubic} Let $S$ be a cubic surface cut out from $\calS_3$ by a hyperplane 
$H: a_0x_0+\cdots+a_5x_5 = 0$. Assume that the intersection is transversal so that $S$ is nonsingular. Then $H$ intersects the planes $\pi_{ab,cd,ef}$ along lines. Thus we have a choice of 15 lines $\ell_{ab,cd,ef}$ from the set of 27 lines on  $S$. The hyperplanes   $H_{ab}$ define 15 tritangent planes $P_{ab}$ of $S$ each containing 3 lines $\ell_{ab,cd,ef}$. Each line is contained in 3 tritangent planes. Thus we get a realization of our $15_3$-configuration as a subconfiguration of the configuration $(27_5,45_3)$ of 27 lines and 45 tritangent planes. The lines $\ell_{ab,cd,ef}$ and their points of intersection form a $(15_6,45_2)$-subconfiguration of the $(27_{10},135_2)$-configuration of lines and points on a cubic surface (one assumes that the cubic surface does not contain Eckardt points).

Two tritangent planes $P_{ij}$ and $P_{mn}$ form a {\it Cremona pair} (see \cite{Reye}, p. 218)  if they do not intersect along one of the 15 lines $\ell_{ab,cd,ef}$, or, equivalently, if $\{i,j\}\cap \{m,n\}\ne \emptyset$. Each plane $P_{ij}$ enters in 8 Cremona pairs. This implies that there are 60 Cremona pairs. The intersection line of two planes from a Cremona pair is called a {\it Pascal line}.  Thus tritangent planes and Pascal lines form a $(15_8,60_2)$-configuration.
It is cut out by $H$ from the dual  plane of a line  joining two diagonal points not on a transversal line.

Each Pascal line intersects the cubic surface $S$ in 3 points. Each point is the intersection point of two lines on $S$ lying in each of the tritangent planes from the pair. Thus through each of 45  intersections points of lines $\ell_{ab,cd,ef}$ passes four Pascal lines. This defines a $(45_4,60_3)$-configuration.  The duals  of the 2-faces $\Pi_I, \#I = 3,$ is a set of 20 lines in $V$. Intersecting with the hyperplane $H$ defines a set of 20 points in $H$ ({\it Steiner points}). The Pascal lines meet by threes in twenty Steiner points. This defines a $(60_1,20_3)$-configuration.

The duals of the edges $\Pi_I, \#I = 2,$ is a set of 15 planes in $V$. Each plane contains 4 duals of the 2-faces. Intersecting with the hyperplane $H$ we find 15 lines in $H$ ({\it Pl\"ucker lines}). Each Pl\"ucker line contains 4 Steiner points.  We get a $(15_4,20_3)$-configuration of lines and points in $\bbP^3$.

The set of 60 Pascal lines is divided into 6 sets of 10 lines. Each set of ten lines is equal to the set $\calP_a$ of lines $H_{ab}\cap H_{ac}$ with $a$ fixed and $c,d\in [1,6]\setminus \{a\}$. The lines from $\calP_a$ are the ten edges of the Sylvester pentahedron of the cubic surface $S$. It forms a $10_3$-subconfiguration of the $(45_4,60_3)$-configuration.  It is isomorphic to the Desargues configuration.

Consider the projection of the cubic surface to the plane from a general point in $\bbP^3$. We get a set $\calL$ of 15 lines (the projections of the 15 lines on the cubic surface), a set $\calP as$ of 60 lines  (the projections of Pascal lines),   a set $\calP lu$ of 15 lines (the projections of Pl\"ucker lines), a set $\calS t$ of 20 points (the projections of 20 Steiner points), a set $\calI$ of 45 points (the projections sof intersection points of the 15 lines on the cubic surface), a set $\calT$ of 15 triangles formed  by the projection of three lines in 15 tritangent planes. 

 The 6 lines from a Cremona pair of tritangents form a hexagon such that its opposite sides intersect at three points lying on a line (the projection of the Pascal line). So this is similar to the Pascal Hexagram, although there is no conic which is inscribed in the hexagon!

The following table summarizes the configurations derived from the Cremona-Richardson configurations in $\bbP^3$:

\bigskip
{\small
\begin{tabular}{ | l | r | r | r | r | r| r| }
\hline
 &$\calL$&$\calP as$&$\calP lu$&$\calS t$&$\calI$&$ \calT$\\
\hline
$\calL$&{}&$(15_{24},60_6)$&&{}&$(15_6,45_2)$&$15_3$\\
\hline
$\calP  as$&$(60_6,15_{24})$&{}&&$(60_1,20_3)$&$(60_3,45_4)$&$(60_2,15_8)$\\
\hline
$\calP lu$&{}&{}&{}&$(15_4,20_3$&&\\
\hline
$\calS t$&{}&$(20_3,60_1)$&$(20_3,15_4)$&&{}&\\
\hline
 $\calI$&$(45_2,15_6)$&$(45_4,60_3)$&{}&{}&&\\
\hline
$\calT$&$15_3$&$(15_8,60_2)$&{}&{}&&\\
\hline
\end{tabular}}

\bigskip

Let set of 15 lines $\ell_{ab,cd,ef}$ is complementary to a double-six of lines. In fact, it is easy to see that no three of the remaining 12 lines can be coplanar, and this characterizes a double-six. 
A general cubic surface can be written in 36 ways as a hyperplane section of the Segre cubic (or equivalently writing its equation in Cremona's hexagonal form:
$$\sum_{i=0}^5l_i^3 = \sum_{i=0}^5l_i = 0,$$ 
where $l_i$ are linear forms in variables $z_0,z_1,z_2,z_3.$)

\subsection{Nodal cubic} Let $H$ be a hyperplane which is tangent to the Segre cubic at one nonsingular point. It cuts out a cubic surface with a node. There are 6 lines passing through the node. The remaining 15 lines are the lines $\ell_{ab,cd,ef}$. Projecting from the node we get a conic with 6 points on it. The blow-up of these six points is a minimal resolution of the cubic. The 15 lines join the six points by pairs. There are 60 hexagons  whose sides belong to this set. Each hexagon defines the Pascal Hexagram.

\subsection{The Segre quartic} Let $\calS_4\subset \bbP_K^4$ be the dual quartic of the Segre cubic (also known  as the {\it Igusa quartic} because of its modular interpretation given by J. Igusa, see \cite{DO}). It contains the transversals as its double lines. The dual of a singular point of the Segre cubic is a hyperplane (called a {\it cardinal space} in \cite{Rich}). They are parametrized by pairs   of complementary  subsets   of $[1,6]$ of cardinality 3.  A cardinal space $C(I)$ corresponding to $I,\bar{I}$  contains 6 transversal lines connecting  9 diagonal points $p_A$, where $A\not\subset I, A\not\subset \bar{I}$. For example, $C(123)$ contains the transversal lines 
$$\ell_{14,25,36},\ \ell_{14,26,35},\ \ell_{15,24,36},\ \ell_{15,26,34},\ \ell_{16,25,34},\ \ell_{16,24 ,35}.$$
The 6 transversals in a cardinal space lie on a unique quadric in this space contained in 
 the Segre quartic  (a {\it cardinal quadric}).  Each transversal line lies in 4 cardinal quadrics. A pair of skew transversal lines lies in a unique quadric, a pair of coplanar transversal lines lie in two quadrics.  Transversal lines and cardinal spaces  define a $(15_4,10_6)$-configuration  which is dual to the $(15_4,10_6)$-configuration of nodes and planes coming from the Segre cubic.

 In the interpretation of $\calS_4$ as a compactification of the moduli space of principally polarized abelian surfaces of level 2, the diagonal points correspond to 0-dimensional boundary components and the  transversal lines correspond to the closures of 1-dimensional boundary components. We  encounter the $15_3$-configuration as the familiar configuration of isotropic planes and lines in the symplectic space $\bbF_2^4$ (see \cite{Geer}). 

We have 20 planes  $\Pi_{I}, \#I = 3$.  Each plane contains 3 lines joining two diagonal points on two edges $\Pi_{J}, J\subset I$ (also called Pascal lines in \cite{Rich}). Through each diagonal point $p_{ab}$ passes three transversals $\ell_{ab,cd,ef},\ell_{ab,ce,df},\ell_{ab,cf,de}$. Each pair spans a plane  which contains 4 Pascal  lines (e.g. the plane spanned by the pair  $\ell_{ab,cd,ef},\ell_{ab,ce,df}$ contains the Pascal lines 
$\langle p_{ce},p_{cd}\rangle, \langle p_{cd},p_{df}\rangle, \langle p_{ef},p_{ce}\rangle,\langle p_{df},p_{ef}\rangle$). The two transversals contained in the plane belong to a cardinal hyperplane. Thus we encounter the $(45_3,60_4)$-configuration of lines and planes which is dual to the $(45_3,60_4)$-configuration considered in \ref{cubic}. 

The  intersection of $\calS_4$ with a transversal hyperplane $H$ is a quartic surface with 15 nodes. There are 10 sets of six coplanar nodes lying on a conic which correspond to cardinal quadrics.    This defines a point-plane $(15_6,10_4)$-configuration. The hyperplane $H$ intersects the 60  Pascal lines at 60 different points ({\it Pascal points}). Also $H$ intersects the 20 planes $\Pi_I,\#I = 3$ at 20 lines (the {\it Steiner lines}). Each Steiner line contains 3 Pascal points.  Each  Pascal  point is contained in 4 Steiner lines. Thus Pascal points and Steiner lines define  a $(60_1,20_3)$-configuration. The edges $\Pi_I,\#I = 2$ define 15 points, the Pl\"ucker points lying by three in 20 Steiner lines. Through each Pl\"ucker point passes 4 Steiner lines.  This defines a $(15_4,20_3)$-configuration of lines and planes in $\bbP^3$. We have 45 lines joining two nodes which contain 4 Pascal points ({\it diagonal lines}). Other 60 pairs of nodes are joined by lines with no Pascal points on it.

 \bigskip
 Let $\calN$ be the set of the 15 nodes of the quartic, $\calC$ be the set of the projections of the  10 conics, $\calP as$ be the set  of the 60 Pascal points, $\calS t$ be the set of the 20 Steiner lines,  $\calP lu$ be the set of  the 15 Pl\"ucker planes and $\calD$ be the set of  the diagonal lines. 

The following table summarizes the configurations derived from the Cremona-Richardson configurations in $\bbP^3$:

\medskip
{\small
\begin{tabular}{ | l | r | r | r | r | r| r|}
\hline
 &$\calN$&$\calC$&$\calP as$&$\calS t$&$\calP lu$&$\calD$\\
\hline
$\calN$&{}&$(15_{4},10_6)$&{}&{}&&$(15_6,45_2)$\\
\hline
$\calC$&$(10_6,15_{4})$&{}&&&&$(10_9,45_2)$\\
\hline
$\calP as$&{}&{}&{}&$(60_1,20_3)$&{}&$(60_3,45_4)$\\
\hline
$\calS t$&{}&&$(20_3,60_1)$&&$(20_3,15_4)$&\\
\hline
 $\calP lu$&&&{}&$(15_4,20_3)$&&\\
\hline
$\calD$&$(45_2,15_6)$&$(45_2,10_9)$&$(45_4,60_3)$&&&\\
\hline
\end{tabular}}

\subsection{Finite geometry interpretation}\label{finite} 
Let $V = \bbF_2^{[1,6]}$ be the linear space  of functions $[1,6]\to \bbF_2$ identified with subsets of $[1,6]$. We have $A+B = A\cup B\setminus A\cap B$. We equip $V = \bbF_2^{[1,6]}$ with a symmetric bilinear form defined by 
\begin{equation}\label{pairing}
\langle A,B\rangle = \# (A\cap B) \mod 2.
\end{equation}
Let $V_0$ be the subspace of subsets of even cardinality. The restriction of the symmetric form to $V_0$ is degenerate. Its radical $R$ is one-dimensional subspace $\{\emptyset,[1,6]\}$. Let $W = V_0/ R$. It is a 4-dimensional vector space over $\bbF_2$ equipped with a non-degenerate symmetric bilinear form $e:W\times W\to \bbF_2$  induced  by  \eqref{pairing}. Since $e(A,A) = 0$ for any subset $A$ of even cardinality, the symmetric form $e$ is a nondegenerate symplectic form. It follows from \eqref{gras} that there are 25 planes in $W$. Fifteen of them are isotropic (i.e. the restriction of $e$ to the subspace is trivial). Each isotropic plane consists of 4 elements represented by subsets $\emptyset, \{a,b\},\{c,d\},\{e,f\}$, where $[1,6]= \{a,b,c,d,e,f\}$. Thus the transversals $\ell_{ab,cd,ef}$ correspond to isotropic planes.  We have $10$ anisotropic (i.e., non-isotropic) planes represented by sets $\{\emptyset,\{a,b\},\{b,c\},\{a,c\}$, where $\{a,b,c\}$ is a subset of 3 elements in $[1,6]$. Of course, replacing $\{a,b,c\}$ by the complementary set does not change the plane. Thus cardinal subspaces correspond to anisotropic planes. Also, observe that non-zero elements which can be represented by 2-elements subsets $I$ of $W$ correspond to diagonal points $p_{I}$. 

We denote an isotropic plane by $L_{ab,cd,ef}$ and an anisotropic plane by $L_{abc}$. It is clear that 
$$L_{ab,cd,ef}\cap L_{abc} =\{0\}\Leftrightarrow \ell_{ab,cd,ef}\subset C(abc).$$
$$I\in L_{ab,cd,ef}\Leftrightarrow p_I\in \ell_{ab,cd,ef}.$$
This gives a finite geometry interpretation of the $(15_4,10_6)$-configuration formed by transversals and cardinal spaces and $15_3$-configuration of diagonal points and transversals. Notice that the symmetry group $S_6$ of the configurations is isomorphic to the symplectic group $\Sp(4,\bbF_2)$ via its natural action on subsets of $[1,6]$.

\section{Kummer configurations}

\subsection {Kummer surface} Intersecting the Segre quartic with a  hyperplane tangent to $\calS_4$ at its nonsingular point defines a quartic with 16 nodes. It is known (\cite{Hudson}) that each 16-nodal quartic surface in $\bbP_K^3$ can be obtained in this way (a modern proof of this classical fact can be found   in \cite{Geer}). It is isomorphic to the Kummer variety of the Jacobian  of a nonsingular curve of genus 2.  Together with 10 conics coming from the cardinal hyperplanes we get 16  conics of the Kummer surface which are cut out by a plane (a {\it trope}). Each conic contains 6 nodes, and each node is contained in 6 tropes. Each pair of nodes is contained in two tropes. This is the famous Kummer $16_6$-configuration. 

The Kummer configuration is a symmetric design of type $\lambda = 2$ (a {\it biplane}). It follows from equation \eqref{design} that for each $v_k$-biplane, $v = 1+\frac{k(k-1)}{2}$. There are three non-isomorphic biplanes with $v = 16$ (see  \cite{Maria},\cite{Hugh}). The Kummer biplane has the largest symmetry group.

 Projecting from the new node, we obtain the familiar configuration of 6 lines in the plane tangent to a conic $C$. The images of the 15 nodes are the intersection points of these lines. The images of the 10 conics are the conics passing through six intersection points no two of which lie on the same line. The six points form the set of vertices of a hexagon which circumscribe $C$. The projection of the 45 diagonals are the diagonals of the hexagons, three in each  hexagon. The Pascal points are the intersection points of three diagonals in the same hexagon which agrees with the Pascal Theorem. 

\subsection{Kummer variety} Let $A$ be an abelian variety of dimension $g$ over an algebraically closed field $K$ of characteristic $\ne 2$. The {\it Kummer variety} of $A$ is defined to be the quotient $\Kum(A) = A/(\iota)$, where $\iota$ is the automorphism $a\mapsto -a$.
It is a normal algebraic variety with $2^{2g}$ singular points corresponding to the fixed points of $\iota$. Each singular point is formally isomorphic to the singular point of the affine cone over the Veronese variety $v_2(\bbP^{g-1})$. 

Assume that $A$ admits an irreducible principal polarization defined by an irreducible divisor $\Theta$ with $h^0(\Theta) = 1$ (e.g. is  $A$ is a Jacobian).  We will also assume that $\Theta$ is symmetric, i.e., $\iota(\Theta) = \Theta$. 
The linear system $|2\Theta|$ maps $A$ to $\bbP^{2^g-1}$ and factors through an embedding  
$$j:\Kum(A)\hookrightarrow  \bbP^{2^g-1}$$
  (\cite{Lange}). The degree of the variety $j(\Kum(A))$ in $\bbP^{2^g-1}$ is equal to $2^{g-1}g!$.

When $g = 2$, we get a Kummer quartic surface in $\bbP^3$ with 16 nodes. 

\subsection{The Kummer $2^{2g}_{2^{g-1}(2^g-1)}$-configurations} For every $2$-torsion point $x\in A$ we can define the translate $\Theta_a = t_a(\Theta)$, where $t_x:a\mapsto a+x$ is the translation automorphism. It is known that $\Theta$ passes through $2^{g-1}(2^g-1)$ $2$-torsion points with odd multiplicity. Let $\calA=\calB$ and $R$ be the relation on $\calA\times \calA$ defined by $(x,y)\in \calA$ if and only if $\textrm{mult}_x\Theta_y $ is odd. This defines a $2^{2g}_{2^{g-1}(2^g-1)}$-configuration. 

Consider an embedding $j:\Kum(A)\hookrightarrow \bbP^{2^g-1}$. The  image of $2$-torsion points is a set of $2^g$ points in the projective space. The image of a divisor $\Theta_a$ is a subvariety of $j(\Kum(A))$ which is cut out by an everywhere tangent hyperplane (a trope). So we have $2^{2g}$ points and hyperplanes in $\bbP^{2^g-1}$ realizing the Kummer $2^{2g}_{2^{g-1}(2^g-1)}$-configuration.

\subsection{Finite geometry interpretation}\label{fin} We generalize the construction from \ref{finite}. Let $W_g$ be the linear space of dimension $2g$ formed by subsets of even cardinality of the set $[1,2g] = \{1,\ldots,2g\}$ with complementary subsets identified. It is  equipped with a nondegenerate symplectic form $e_g$ defined by \ref{pairing}. Let $Q_g$ be set of subsets of $[1,2g]$ of cardinality $\equiv g+1\mod 2$ with complementary subsets identified.
It has a natural structure of an affine space  over $W_g$ with respect to the symmetric sum of subsets. We will identify $Q_g$ with quadratic forms on $W_g$ by setting, for any $S\in Q_g$, and any $T\in W_g$,
$$q_S(T) = \frac{1}{2}\#T+\#T\cap S \mod 2 = \frac{1}{2}\#T+e_g(S,T) \mod 2.$$
The associated symmetric bilinear form of $q_S$ coincides with $e_g$.  Consider the relation $R\subset W_g\times Q_g$ defined by
$$(T,S)\in R\Leftrightarrow q_S(T) =
\begin{cases}
      1& \text{if}\  \#S\equiv g+1\mod 4, \\
      0& \text{if}\ \#S\equiv g-1\mod 4.
\end{cases} 
$$
Recall that there are two types of nondegenerate quadratic forms on $(\bbF_2)^{2g}$. They are distinguished by the number of isotropic vectors $x\in \bbF_2^{2g}$. For the {\it even type} this number is equal to $2^{g-1}(2^g+1)$ and for the {\it odd type} it is equal to 
$2^{g-1}(2^g-1)$. One checks that the quadratic form $q_S$ is of even type if and only if $\#S\equiv g+1\mod 4$. This implies that $\#R(S) = 2^{g-1}(2^g-1)$ for any $S\in Q_g$.
The relation $R$ has the obvious symmetry group isomorphic to $\Sp(2g,\bbF_2)$. This implies that $\#R(T) = 2^{g-1}(2^g-1)$ for any $T\in W_g$. Thus we get a $2^{2g}_{2^{g-1}(2^g-1)}$-configuration with symmetry group $\Sp(2g,\bbF_2)$. There is more symmetry. Let us see that translations by elements from $W_g$ are symmetries of the configuration. For any $T\in W_g$, we have
\[q_{S+T}(T') = q_S(T')+e_g(S+T,T')+e_g(S,T') = q_S(T')+e_g(T,T').\]
Thus
\[q_{S+T}(T+T') = q_S(T+T')+e_g(T,T+T') = \]
\[q_S(T)+q_S(T')+e_g(T,T')+e_g(T,T+T') = 
q_S(T)+q_S(T').\] 
Obviously,  $\#(S+T)-\#S =\#T -2\#(S\cap T)$. 
Suppose $\#T\equiv 0 \mod 4$. Then 
$\#(S+T)-\#S =2\#(S\cap T)\mod 4$ and $q_S(T) = \#(S\cap T)\mod 2$.
This implies that $(T',S)\in R$ if and only if $(T+T',S+T)\in R$. 

Suppose $\#T\equiv 2 \mod 4$. Then 
$\#(S+T)-\#S =2+2\#(S\cap T)\mod 4$ and $q_S(T) = 1+\#(S\cap T)\mod 2$.
Again this implies that $(T',S)\in R$ if and only if $(T+T',S+T)\in R$. 
This checks that $W_g$ is a subgroup of the group of symmetries. Thus, we found that 
the semi-direct product $G = \bbF_2^{2g}\rtimes \Sp(2g,\bbF_2)$ is a group of symmetries.

It is proven by A. Krazer \cite{Krazer}  (for other expositions of this result  see \cite{Coble}, \cite{DoOr}, \cite{Mumford}) that the obtained configuration is isomorphic to the Kummer configuration.

\subsection{Realizations of symmetries}\label{schro} The group of symmetries of the Kummer configuration defined by translations can be realized by projective transformations of the space $\bbP^{2^g-1}$.  Recall the Schr\"odinger representation of the group $\bbF_p^2$ in $\bbP^{n-1}$ which we discussed in \ref{segre}. This can be generalized to a projective representation of $\bbF_p^{2g}$ in $\bbP^{p^g-1}$. First we define a central extension $\calH_g(p)$ of the group $\bbF_p^{2g}$ with the center isomorphic to the group $\mu_p$ of complex $p$th roots of unity  (the {\it Heisenberg group}). It is defined by the following group law on the set $\bbF_p^{2g}\times \bbF_p$: 
\[(v,\alpha)\cdot (w,\beta) = (v+w,e^{2\pi \langle v,w\rangle}\alpha\beta),\]
where $\langle {},{}\rangle:\bbF_p^{2g}\times \bbF_p^{2g}\to \bbF_p$ is the standard symplectic form defined by the matrix $\begin{pmatrix}0_g&I_g\\
-I_g&0_g\end{pmatrix}$.
The center $C$ is the subgroup generated by $(0,1)$, the quotient group is isomorphic to 
$\bbF_p^{2g}$. Let $(e_1,\ldots,e_{2g})$ be the standard basis in $\bbF_p^{2g}$. Let $A$ be the subspace spanned by the first $g$ unit vectors, and $B$ be the subspace spanned by the last $g$ unit vectors. These are maximal isotropic subspaces. Let $\bar{A} = \{(v,1),v\in A\}, \bar{B} = \{(w,1),w\in B\}$. These are subgroups of $\calH_g(p)$ isomorphic to $\bbF_p^{g}$. For any $(v,1)\in \bar{A}$ and $(w,1)\in \bar{B}$, we have
$$[(v,1),(w,1)] = (0,e^{2\pi i \langle v,w\rangle}).$$
Let us write any element in $\calH_g(p)$ as a triple $(x,y,\alpha)$, where $x\in \bar{A}, y\in \bar{B}$ and $\alpha\in C$. The Schr\"odinger linear representation of $\calH_g(p)$ in $k^{p^g}$ is defined by the formula
\[\bigl((x,y,\alpha)\cdot \phi\bigr) (z)= \alpha e^{2\pi i\langle z,y\rangle}\phi(z+x),\]
where we identify the linear space  $K^{p^g}$ with the space of $K$-valued functions on the set $A = \bbF_p^g$. The Schr\"odinger projective representation in $\bbP^{p^g-1}$ is the corresponding representation of $\calH_g(p)$ in $\bbP(K^{p^g})$. It is clear that the center acts trivially in the projective space, so we have a projective representation of the group $\bbF_p^{2g}$ in $\bbP^{p^g-1}$.

Assume $p = 2$. The group of $2$-torsion points of an abelian variety $A$ of dimension $g$ is isomorphic to $\bbF_2^{2g}$. It is also equipped with a canonical symplectic form defined by the Weil pairing. This allows to define its center extension isomorphic to the Heisenberg group $\calH_g(2)$ (the {\it theta group}). If $\Theta$ is a divisor on $A$ defining a principal polarization, then the theta group acts linearly in the space $L(2\Theta) = H^0(A,\calO_A(2\Theta))$ via its natural action on $A$ by translations. This defines a $\calH_g(2)$-equivariant embedding of the Kummer variety in $\bbP^{2^g-1}$ (see \cite{Lange}, Chapter 6). Thus the group $\bbF_2^{2g}$ acts on the Kummer variety and hence is realized as a group of symmetries of the Kummer configuration. One can show that this action is isomorphic to the action of $W_g$ on the Kummer configuration defined in \ref{fin}.

\subsection{The Kummer designs} The Kummer configuration  is a symmetric design of type $\lambda = 2^{g-1}(2^{g-1}-1)$. Let us prove this by using the construction from the previous section. 
Let $T,T'$ be a pair of  elements of $W_g$. Since $\Sp(2g,\bbF_2)$ is 2-transitive on the set of pair of nonzero vectors $A,B$ with given value of $e_g(A,B)$, it suffices  to consider the following two cases: 

Case 1: $T = \{1,2\}, \ T' = \{2,3\}, \ e_g(T,T') = 1$

Case 2: $T = \{1,2\},\  T' = \{3,4\}, \ e_g(T,T') = 0$.

Let $\lambda_1 = \#R(T)\cap R(T')$ in the first case and $\lambda_2 = \#R(T)\cap R(T')$ in the second case. By symmetry, it is enough to show that $\lambda_1 = \lambda_2$, and then $\lambda$ could be found from equation \eqref{design}. 

Consider the first case. Let us count the number of sets $S\in Q_g$ such that $T,T'\in R(S)$. Suppose 
$\#S \equiv g+1\mod 4$. Then $T,T'\in R(S)$ if and only if 
$\#T\cap S = \#T'\cap S = 0 \mod 2.$ This is equivalent to that $S$ or $\bar{S}$  does not contain $1,2,3$. Denote the set of such subsets by $\calP_1$. Now suppose that $\#S \equiv g-1\mod 4$. Then 
$T,T'\in R(S)$ if and only if $\#T\cap S = \#T'\cap S = 1 \mod 2$. This is equivalent to that $S$ or $\bar{S}$ contains $2$ and does not contain $1,3$. Denote the set of such subsets by $\calP_2$

Consider the second case. Suppose 
$\#S \equiv g+1\mod 4$. Then $T,T'\in R(S)$ if and only if 
$\#T\cap S = \#T'\cap S = 0 \mod 2.$ This is equivalent to that $S$ or $\bar{S}$  does not contain $1,2,3,4$ or does not contain $\{1,2\}$ but contains $\{3,4\}$. Denote the set of such subsets by $\calQ_1$. Now suppose that $\#S \equiv g-1\mod 4$. Then 
$T,T'\in R(S)$ if and only if  $S$ or $\bar{S}$ does not contain 1, contains  $2$ and one of the numbers $3,4$.  Denote the set of such subsets by $Q_2$. 

Let $S\in \calP_1$. We may assume that $S\subset [4,\ldots,2g+2]$. If $4\not\in S$, then 
$S\in Q_1$. If $4\in S$, then adding $2$ to $S$ and  one of the numbers from $[5,\ldots,2g+2]$  we get an element in $Q_2$. Let $S\in \calP_2$. We may assume that $2\in \calP_2$. If $S$ contains $4$, then $S\in \calQ_2$. If $4\not\in S$, then deleting  4  and one of the numbers from $[5,\ldots,2g+2]$ we obtain an element from $\calQ_1$. This shows that 
$$\lambda_1 = \#\calP_1+\#\calP_2 \le \lambda_2 =\#\calQ_1+\#\calQ_2.$$
Similarly, we prove that $\lambda_2\le \lambda_1$.

\subsection{Hadamard designs} The Kummer design is a {\it Hadamard design}. Recall that an {\it Hadamard matrix} is a square matrix $A$ of order $n$ with $\pm 1$  whose rows  are orthogonal. 
It is easy to see that $n=4t$ is a multiple of 4. By permuting rows and columns, one may assume that the first column and the first row does not contain $-1$. Consider the matrix $A'$ obtained from $A$ by deleting the first row and the first column. We take $\calA$ to be the set of rows of $A'$ and $\calB$ the set of columns of $A'$. We define the incidence relation by requiring that a row $A_i$ is incident to the column $B_j$ if the entry $a_{ij}$ of $A'$ is equal to 1. This defined a symmetric $4t-1_{2t-1}$-design of type $\lambda = t-1$. Any design of this type is obtained from an Hadamard matrix. An example of an Hadamard design is the Fano plane $(t = 2)$.  The Kummer design is not an Hadamard design however its $(-1,1)$-incidence matrix is an Hadamard matrix with $t = 2^{g-2}$ (see \cite{Kantor2}, p.47).

\subsection{2-transitive designs} Note that according to a result of W. Kantor \cite{Kant},  there are only the following symmetric designs (up to taking the complementary design) whose symmetry group acts $2$-transitively on the set of points (or blocks):

1) A point-hyperplane design in $\PG(n,q)$;

2) The unique Hadamard $11_5$-design with $t = 3$. Its proper symmetry group is isomorphic to $\PGL(2,\bbF_{11})$ (see \cite{Kantor2}, p. 68).

3) A certain $176_{50}$ design with $\lambda = 14$ constructed by G. Higman (see \cite{Higman}). Its proper symmetry group is isomorphic to the Higman-Sims simple group of order 
$2^93^25^37.11$.

4) The Kummer designs.

It is an interesting problem to realize designs 2) and 3) in algebraic geometry.

 \subsection{Klein's $60_{15}$-configuration} (see \cite{BN}, \cite{Hudson}) Consider the Schr\"odinger linear representation of the group $\calH_2(2)$ in $K^4$ defined in \ref{schro}. We choose the coordinates in $K^4$ such that the unit vectors $e_1,e_2,e_3,e_4$ correspond to the characteristic functions of $\{0\},\{(1,0)\},\{(0,1)\},\{(1,1)\}$, respectively. The generators of the subgroups $\bar{A}$ and $\bar{B}$ act by the formula
\begin{eqnarray*}
(e_1,1)\cdot (z_0,z_1,z_2,z_3) &= &(z_1,z_0,z_3,z_2),\\
(e_2,1)\cdot (z_0,z_1,z_2,z_3) &= &(z_2,z_3,z_0,z_1),\\
(e_3,1)\cdot (z_0,z_1,z_2,z_3) &= &(z_0,-z_1,z_3,-z_2),\\
(e_4,1)\cdot (z_0,z_1,z_2,z_3) &= &(z_0,z_1,-z_3,-z_2).
\end{eqnarray*}
For any $g\in\calH_2(2)\setminus\{1\}$ its set of fixed points is a pair of  planes $L_g^{\pm}$ which define a pair of skew lines $l_g^\pm$ in $\bbP^3$. The decomposition $K^4 = L_g^++ L_g^-$ is the decomposition into the direct sum of eigensubpaces of $g$. If $g^2 = 1$, the eigenvalues are $\pm 1$. If $g^2 = -1$, the eigenvalues are $\pm i$. Let $L$ be any plane in $\bbF_2^{2g}$ and $\bar{L}$ be its preimage in $\calH_2(2)$. If $L$ is isotropic with respect to $\langle{},{}\rangle$, then all elements in $\bar{L}$ commute. Thus, if $g,g'\in \bar{L}, g\ne \pm g'$, the sets of fixed points of $g$ and $g'$ have common points, two lines in each plane. This shows that three pairs of skew lines $l_g^\pm$ defined by the  nonzero elements from $L$ form a tetrahedron $T_L$. The opposite edges of the tetrahedron are the fixed lines of a nonzero element from $L$. 

Two isotropic planes may have either one common nonzero vector in common (an {\it azygetic pair}) or no common nonzero vectors (a {\it syzygetic pair}). The invariant tetrahedra of an azygetic pair have a pair of opposite edges in common. The invariant tetrahedra of a syzygetic pair have no edges in common. For every isotropic plane $L$, there are 3 isotropic pairs $L'$ which have a fixed common nonzero vector $v\in L$. This shows that each edge is contained in 3 tetrahedra, so we have 30 edges. Also we see that no two tetrahedra share a vertex. Thus we have a set $\calA$ of 60 vertices and a set $\calB$ of 60 faces of the 15 tetrahedra. Six faces pass through an edge (two from each of the three tetrahedra containing it). Each face contains 3 edges. This defines a $(60_3,30_6)$-configuration. Each edge contains 6 vertices (two from each tetrahedra it contains). Since a face contains 3 edges, each has a common vertex, we have 15 vertices in each face. One can also see that each vertex is contained in  15 faces. Thus faces and vertices form a $60_{15}$-configuration. This is Klein's $60_{15}$.

If $L$ is anisotropic, the sets of fixed points of nonzero elements have no common nonzero elements. Thus the corresponding three pairs of lines are skew. An orthogonal anisotropic plane $L^\perp$  defines another 6 skew lines such that each fix-line defined  by $L$  intersects each fix-line defined by $L^\perp$. Let $Q$ be the unique quadric in $\bbP^3$ which contains 3 fix-lines defined by $L$. Then each fix-line defined by $L^\perp$ must lie in $Q$. Then the remaining 3 fix-lines defined by $L$ also must lie in $Q$. Thus $Q$ contains 12 fix-lines defined by the pair $(L,L^\perp)$. There are 10 such pairs, so we obtain 10 quadrics. They are called {\it  Klein's fundamental quadrics}.

Consider the tensor square of the Schr\"odinger linear representation in the dual space $(K^4)^*$. It decomposes into the direct sum $\Lambda^2(K^4)^*\oplus S^2((K^4)^*)$. The center of the Heisenberg group acts identically. So, we get a linear representation of an abelian group $\bbF_2^4$. It decomposes into the direct sum of one-dimensional representations. The 10 fundamental quadrics represent a basis of eigensubspaces in $S^2((K^4)^*)$. Let $\omega_i, 1\le i\le 6,$ be a basis of eigensubspaces in $\Lambda^2(K^4)^*$ (the {\it six fundamental complexes}). Each defines a polarity, i.e. a linear isomorphism $\omega_i:K^4\to (K^4)^*$. The Kummer configuration can be obtained as follows. Take an arbitrary point $P\in \bbP^3$. The orbit of $P$ with respect to the Schr\"odinger representation gives a set of 16 points. The images of $P$ under $\omega_i$ define 6 planes. The polar plane of $P$ with respect to the fundamental quadrics (defined by the orthogonal subspace of the corresponding line in $K^4$) define another 10 planes. Thus we get a set of 16 points and a set of 16 planes. There exists a unique Kummer surface in $\bbP^3$ which realizes these sets as the set of 16 nodes and 16 tropes.

Recall that each element in $\Lambda^2((K^4)^*)$ defines a hyperplane section of the Grassmannian $G_{1,3}$ in its Pl\"ucker embedding. Thus each subset of 4 fundamental complexes $\omega_i$ define a pair of lines in $\bbP^3$. Together we get 15 pairs of lines. This is the set of 30 edges of the 15 fundamental tetrahedra. In this way the Klein $60_{15}$-configuration can be reconstructed from the six fundamental compexes.

\subsection{A Kummer $16_{10}$-configuration} A Jacobian Kummer surface is obtained as the image of the linear system $|2\Theta|$ on the Jacobian variety of a genus 2 curve. A minimal nonsingular model of the Kummer surface of a simple abelian surface with polarization of type $(1,3)$ admits an embedding in $\bbP^3$ as a quartic surface with two sets of 16 disjoint lines forming a symmetric $16_{10}$-design of type $\lambda = 5$. This surface was first discovered by M. Traynard in 1907 (\cite{Tray}) and rediscovered almost a century later in  \cite{BN}, \cite{Naruki}.  The Kummer $16_6$ and $16_{10}$ configurations are complementary to each other.

\section{A symmetric realization of $\bbP^2(\bbF_q)$}

\subsection{Mukai's realization} We know already  that the projective plane $\bbP^2(\bbF_q)$ can be realized by points and lines in $\bbP_{\bar{\bbF}_q}^2$, or by its blow-up. Although all  these realizations admit the symmetry group $\PGL(3,\bbF_q)$, none of them admits a switch. The following realization admitting a switch is due to  S. Mukai.  Let $\bbP^2 = \bbP_{\bar{\bbF}_q}^2$ and $X_q$ be a surface in $\bbP^2\times \bbP^2$ given by the equations
$$x_0y_0^q+x_1y_1^q+x_2y_2^q = 0, \quad  y_0x_0^q+y_1x_1^q+y_2x_2^q = 0.$$
It is easy to verify that $X_q$ is a nonsingular minimal surface with 
$$q(X_q) = 0, \quad p_g(X_q) = \frac{1}{4}q^2(q-1)^2, \quad K_{X_q}^2 = 2(q-2)^2(q^2+1).$$
Let $p_i:X_q\to \bbP^2, i = 1,2,$ be the two projections. Let $a = (a_0,a_1,a_2)\in \bbP^2(\bbF_{q^2})$. Since 
$$(y_0a_0^q+y_1a_1^q+y_2a_2^q)^q = y_0^qa_0^{q^2}+y_1^qa_1^{q^2}+y_2^qa_2^{q^2} = y_0^qa_0+y_1^qa_1+y_2^qa_2,$$ we see that the fibre $R_a = p_1^{-1}(a)$ is isomorphic to a line  in $\bbP^2$ given by the equation  $a_0^qT_0+a_1^qT_1+a_2^qT_2 = 0.$
  Thus we obtain a disjoint set $\calA$ of $q^4+q^2+1$ smooth rational curves $R_a$ in $X_q$. By the adjunction formula, $R_i^2 = -q$. Similarly, considering the second projection we obtain another set $\calB$ of such curves $Q_b$. A curve $R_a$ intersects $Q_b$ if and only if  $b_0a_0^q+b_1a_1^q+b_2a_2^q = 0$, or, equivalently, 
$a_0b_0^q+a_1b_1^q+a_2b_2^q = 0$. Let $q = p^k$ and let $\bfF:(x_0,x_1,x_2)\to (x_0^p,x_1^p,x_2^p) $ be the Frobenius endomorphism of $\bbP_{\bbF_p}^2$. We see that, for any $a,b\in \bbP^2(\bbF_{q^2}) $, 
$R_a\cap Q_b\ne\emptyset $ if and only if the point $\bfF^k(a)$ lies on the line $b_0T_0+b_1T_1+b_2T_2 = 0$. This shows that the configuration $\calA\cup\calB$ is isomorphic to the configuration $\PG(2,\bbF_{q^2})$ under the  map $\bfF^k$. Of course, it contains the subconfigurations $\PG(2,\bbF_{p^i})$ for any $1\le i\le k$.
 
The group 
$\PGL(3,\bbF_{q^2})$ acting diagonally on $\bbP^2\times \bbP^2$ together with a switch defined by the interchanging the factors of $\bbP^2\times \bbP^2$ realizes the subgroup of index $2k$, where $q = p^k$,  of the group of abstract symmetries of the configuration $\bbP^2(\bbF_q)$. The  cosets of non-realizable symmetries are the cosets of the powers of the Frobenius map $\bfF$.

\subsection{The surface $X_2$} In this case the surface $X_2$ is a K3 surface and the configuration $\PG(2,\bbF_4)$ is realized by two sets of 21 disjoint $(-2)$-curves. This surface is the subject of my joint  work with S. Kondo (\cite{DK}). We prove that the following properties characterize a K3 surface $X_2$ over an algebraically closed  field of characteristic 2:

\begin{itemize}
\item [(i)] $X$ is isomorphic to $X_2$;
\item[(ii)]  The Picard lattice of $X$ is isomorphic to $U\perp
D_{20}
\ ;$
\item[(iii)] $X$ has a jacobian quasi-elliptic fibration with one fiber of
type $\tilde{D}_{20} \ ;$
\item[(iv)] $X$ has a quasi-elliptic fibration with the Weierstrass
equation \[y^2 = x^3+t^2x+t^{11};\]
\item[(v)] $X$ has a quasi-elliptic fibration with 5 fibers of type
$\tilde{D}_4$ and the group of sections isomorphic to $(\bbZ/2)^4\ ;$
\item[(vi)] $X$ contains a set $\calA$ of \ $21$ disjoint $(-2)$-curves and
another set $\calB$ of \ $21$ disjoint
$(-2)$-curves such
that each curve from one set intersects exactly 5 curves from the other set
with multiplicity $1 ;$
\item[(vii)] $X$ is isomorphic to a minimal nonsingular model of  a hypersurface of degree 6 in $\bbP(1,1,1,3)$ given by the equation
of  $\bbP^2$ with branch divisor
\[x_3^2+x_0x_1x_2(x_0^3+x_1^3+x_2^3) = 0;\]
\item[(viii)] $X$ is isomorphic to a minimal nonsingular model of the
quartic surface with $7$ rational double points of type $A_3$ which is
defined by the equation
\[x_0^4+x_1^4+x_2^4+x_3^4+x_0^2x_1^2+x_0^2x_2^2+x_1^2x_2^2+x_0x_1x_2(x_0+x_1
+x_2) = 0.\]
\end{itemize}

\subsection{Automorphism group of $X_2$} The full automorphism group $\Aut(X_2)$ is infinite. It contains a normal infinite subgroup
generated by $168$ involutions and the quotient is a finite group
isomorphic to $\PGL(3,\bbF_4)\cdot 2$.

\subsection{The Ceva(3) again} Consider the subset $S$ of 9 points in $\bbP^2(\bbF_4)$ with all coordinates being nonzero. In Mukai's realization they define a subconfiguration of type $9_3$. The group of automorphisms of the surface $X_2$ generated by a switch and the subgroup of $\PGL(3,\bbF_4)$ which leaves the set $S$ invariant is of order 108. It realizes a subgroup of index 2 of the symmetry group of Ceva(3). The existence of the Ceva configuration in $\bbP^2(\bbF_4)$ was first observed by Gino Fano (\cite{Fano2}).

\end{document}